\documentclass[a4paper,11pt,twoside,reqno]{amsart}
\usepackage{latexsym,amsmath,amsthm,verbatim,ifthen,amssymb,graphicx,color,mathrsfs}
\usepackage[all]{xy}
\usepackage{color}
 \title[Homotopy theory of complete Lie algebras]{Homotopy theory of complete Lie algebras and Lie models of simplicial sets}

\author{Urtzi Buijs}
\address{Departamento de Algebra, Geometr\'ia y Topolog\'ia\\
         Universidad de M\'alaga\\
        Ap. 59\\
         29080-M\'alaga,\\
         Espa\~na}
\email{ubuijs@uma.es}

\author{Yves F\'elix}
\address{Institut de Math\'ematiques et Physique\\
         Universit\'e Catholique de Louvain-la-Neuve\\
         Louvain-la-Neuve\\
         Belgique}
\email{Yves.felix@uclouvain.be}

\author{Aniceto Murillo}
\address{Departamento de Algebra, Geometr\'ia y Topolog\'ia\\
         Universidad de M\'alaga\\
        Ap. 59\\
         29080-M\'alaga,\\
         Espa\~na}
         \email{aniceto@uma.es}

\author{Daniel Tanr\'e}
\address{D\'epartement de Math\'ematiques, UMR 8524\\
         Universit\'e de Lille~1\\
         59655 Villeneuve d'Ascq Cedex\\
         France}
\email{Daniel.Tanre@univ-lille1.fr}

\thanks{The first author has been partially supported by the Ram\'on y Cajal MINECO programme.
The first and third authors have been partially supported by the Junta de Andaluc\'\i a grant FQM-213.
The fourth author has been partially supported by the  ANR-11-LABX-0007-01  ``CEMPI''.
The authors are partially supported by the MINECO grant MTM2013-{41768-P}.}



\SelectTips{cm}{}
\newtheorem{theorem}{Theorem}[section]

\newtheorem{corollary}[theorem]{Corollary}
 \newtheorem{lemma}[theorem]{Lemma}
 \newtheorem{proposition}[theorem]{Proposition}
 \theoremstyle{definition}
 \newtheorem{definition}[theorem]{Definition}
 \theoremstyle{remark}
 \newtheorem{remark}[theorem]{Remark}
 \newtheorem{example}[theorem]{Example}

\def\cA{{\mathcal A}}

\def\cF{{\mathcal F}}
\def\cG{{\mathcal G}}

\def\cI{{\mathcal I}}

\newcommand{\bz}{\mathbb Z}
\newcommand{\bq}{\mathbb Q}

\def\L{\mathbb{L}}

\def\Q{\mathbb{Q}}

\def\Z{\mathbb{Z}}

\def\hL{{\widehat{\mathbb L}}}

\def\ad{{\rm ad}}
\def\id{{\rm id}}

    \newcommand{\lasu}{{\mathfrak{L}}}

    \newcommand{\ov}{\overline}

 \newcommand{\lib }{\mathbb{L}}

 \newcommand{\cyl}{\operatorname{\textrm{Cyl}\,}}

\newcommand{\catcdga}{\operatorname{{\bf CDGA}}}

\newcommand{\ho}{\operatorname{{\rm Ho}}}

\newcommand{\catcdgl}{\operatorname{{\bf cDGL}}}

\newcommand{\catss}{\operatorname{{\bf SSet}}}

\newcommand{\Hom}{\operatorname{\text{\rm Hom}}}

 \newcommand{\coch}{ \mathscr{C}}

 \newcommand{\MC}{\operatorname{{\rm MC}}}
  \newcommand{\cil}{\operatorname{{\rm Cyl}}}
\newcommand{\mc}{{\MC}}

   \newcommand{\libc}{{\widehat\lib}}

\begin{document}

\date{\today}

\begin{abstract}
In a previous work,   by extending the classical Quillen construction to the non-simply connected case, we have built a pair of adjoint functors, {\em model} and {\em realization},
$$\xymatrix{\catss&& {\bf cDGL} \ar@<1ex>[ll]^(.46){\langle-\rangle\mbox{\it (realization)}}
\ar@<1ex>[ll];[]^(.50){\lasu   \mbox{\it (model)}}\\}
$$
between the categories of simplicial sets and complete differential graded Lie algebras.  This paper is a follow up of this work. We show that when $X$ is a finite connected simplicial set, then $\langle \lasu_X\rangle$ coincides with $\bq_\infty X^+$, the disjoint union of the Bousfield-Kan completion of $X$ with an external point. We also define a model category structure on {\bf cDGL} making  $\lasu$ and $\langle-\rangle$ a Quillen pair, and we construct an explicit cylinder. In particular, these functors preserve homotopies and weak equivalences and therefore, this gives the basis for developing a {\bf cDGL} rational homotopy theory  for all spaces.
\end{abstract}

\maketitle

\section*{Introduction}

In \cite{pri}, we have built a pair of adjoint functors, {\em model} and {\em realization},
$$\xymatrix{ \catss& \catcdgl \ar@<1ex>[l]^(.50){\langle-\rangle}
\ar@<1ex>[l];[]^(.50){\lasu}\\}
$$
between the usual category of simplicial sets and the category of complete differential graded Lie algebras (cdgl's) 
defined over $\Q$. This category and the main properties of these functors are described in detail in 
Sections~\ref{sec:cdgl} and \ref{sec:modreal}
respectively. Here we briefly recall that, when $X$ is a finite simplicial complex, 
or more generally, a finite simplicial set,   $\lasu_X$ is naturally isomorphic to a free cdgl $(\libc(V),d)$
in which $V$ together with the linear part of $d$  is the
desuspension of the chain complex of $X$, and any vertex of $X$ corresponds to a Maurer-Cartan
element $a\in V$. In particular, when applied to the simplices $\Delta^n$, the model functor gives a cosimplicial  cdgl
 $\lasu_{\Delta^\bullet}$, which in turn,  defines  the \emph{realization} functor for any cdgl $L$ as the simplicial set
 $\langle L\rangle = \Hom_{\bf cDGL} (\lasu_{\Delta^\bullet}, L)$.

 When $X$ is a finite simply connected simplicial set, the component $\lasu_X^a$ of the model $\lasu_X$ 
 at a   non-trivial Maurer-Cartan element $a$ (see Definition~\ref{def:dglmc}) is quasi-isomorphic 
 to the Quillen rational model of $X$ \cite[Theorem C]{pri}.  
 In particular, in view of equality (\ref{picero}) and Theorems \ref{componentecontractil} and \ref{componentecompleccion},
  we deduce that  $\langle \lasu_X\rangle\simeq X_\bq^+$, the disjoint union 
  of the rationalization of $X$ with an external point. 
  Note that our realization  coincides with the classical Quillen realization in the finite simply connected case 
 \cite[Theorem D]{pri}.

 Our first result extends this to the  non-simply connected case.

 \begin{theorem}\label{realization}
The realization $\langle \lasu_X\rangle$ of the model of a connected finite simplicial set $X$ has the homotopy type of $\bq_\infty X^+$, the disjoint union of the Bousfield-Kan $\bq$-completion of $X$ with an external point.
\end{theorem}

We then use  the transfer principle \cite{ban,bermoer} to endow {\bf cDGL} with a cofibrantly generated  model category structure arising from the one in $\catss$.

\begin{theorem}\label{model}
 There is a cofibrantly generated model category structure on $\catcdgl$ for which:
\begin{enumerate}
\item[$\bullet$]
A morphism $f\colon A\to B$ is a fibration if 
it is surjective in non-negative degrees.
\item[$\bullet$]
A morphism $f\colon A\to B$ is a weak equivalence if 
$\widetilde\mc(f)\colon \widetilde\mc(A)\stackrel{\cong}{\to}\widetilde\mc(B)$ is a bijection and  
$f^{a}\colon A^{a}\stackrel{\simeq}{\to} B^{f(a)}$ 
is a quasi-isomorphism for every 
gauge equivalence class $a\in\widetilde\mc(A)$.
\item[$\bullet$]
 A morphism is a cofibration if it has the left lifting property with respect to trivial fibrations.
 \item[$\bullet$] The cdgl morphisms
$$
\lasu_{\mathscr I}=\{\lasu_{\partial\Delta^n}\hookrightarrow \lasu_{\Delta^n}\}_{n\ge0},\quad
\lasu_{\mathscr J}=\{\lasu_{\wedge^n_i}{\hookrightarrow} \lasu_{\Delta^n}\}_{n\ge0, \,i=0,\dots,n},
$$
are  generating sets of cofibrations and trivial cofibrations respectively.
\end{enumerate}
 \end{theorem}

 As a follow up, we obtain the compatibility of the functors realization and model with the relation of homotopy.

 \begin{theorem} The realization and model functors,
$$
\xymatrix{  \catss \ar@<0.75ex>[r]^-{\lasu} &\catcdgl, \ar@<0.75ex>[l]^(0.54){\langle-\rangle} }
$$
form a Quillen pair. In particular, they induce adjoint functors in the homotopy categories,
$$
\xymatrix{  \ho\catss \ar@<0.75ex>[r]^-{\lasu} &\ho\catcdgl, \ar@<0.75ex>[l]^(0.52){\langle-\rangle} }
$$
and both preserve weak equivalences and homotopies.
\hfill$\square$
\end{theorem}

To go further in the model structure of {\bf cDGL} we construct an explicit cylinder object   
as   a completed version of the ``cylinder" introduced by D.~Tanr\'e
  in \cite[II.5]{tan}.   
   In particular, we characterize the {\em Lawrence-Sullivan interval} as a cylinder. 
   Recall \cite{lawsu}  that the Lawrence-Sullivan interval is the cdgl
  $$ \lasu_{\Delta^1}=(\libc(a,b,x),d)\,$$
in which $a$ and $b$ are Maurer-Cartan elements, $x$ is  of degree $0$, and
$$
 d x = \ad_xb + \sum_{n=0}^\infty \frac{B_n}{n!}\,  \ad_x^n (b-a)
 =\ad_xb + \frac{\ad_x}{e^{\ad_x}-1} (b-a)\,,
$$
 where  $B_n$ is the $n$th   Bernoulli number.

 \begin{theorem}\label{thmm:LST}
 $\lasu_{\Delta^1}$ is isomorphic to  the cylinder $\cil(\lasu_{\Delta^0})$.
 \end{theorem}

 We also present a large class of cofibrations in this model structure which we call {\em free extensions} and we deduce the following.

\begin{proposition}
 For any finite simplicial set $X$, its model $\lasu_X$ is a cofibrant cdgl.
 \end{proposition}
 
 In Section~\ref{sec:modelsset}, we also introduce Lie minimal models of 
 connected simplicial sets and detail some examples.

\medskip
Apart from the usual model structure on categories of unbounded chain complexes enriched 
with a suplementary algebraic structure,
A. Lazarev and M. Markl have defined in \cite{lamar} another  structure of model category on 
a subcategory  of {\bf cDGL}.  We   compare them with our model structure in the last section of the text, proving in
Theorem~\ref{propcompar} that the two notions of weak equivalences
 coincide for a large class of maps, including the models of simplicial maps.

\medskip
 The purpose of this paper is to develop the basis of a rational homotopy theory for all spaces. 
 For simply connected spaces this was done in the seminal papers of Quillen \cite{qui} and Sullivan \cite{su}, 
 see also \cite{bousgu}, \cite{FHTI}.
  Among  the main results there are equivalences  between the homotopy category of simply connected spaces 
  with finite Betti numbers and on one side a homotopy category of commutative differential graded algebras, 
  and on the other side a homotopy category of  differential graded Lie algebras. 
  Both approaches have also been extended to spaces with finitely connected components, 
  each of them being a nilpotent space with finite Betti numbers, and the equivalence of homotopy categories 
  extended to this context \cite{lamar}. See also \cite{BM2} for a global theory of  non-connected spaces.
  This work is a step forward to extend the Quillen approach to all spaces.

  \medskip
Mention that Sullivan theory of minimal models has been very efficient for proving new results 
in algebraic topology and differential geometry. 
 Let us recall for instance the dichotomy theorem for a simply connected finite CW-complex $X$: 
 either $\pi_n(X)\otimes \mathbb Q= 0$ for $n$ enough large, or 
 else the sequence $\sum_{q\leq n}\mbox{rank}\, \pi_q(X)$ has an exponential growth (\cite[Theorem 33.2]{FHTI}). 
 In geometry, one of the first application is the Vigu\'e-Sullivan Theorem (\cite{MR0455028}): 
 If $M$ is a simply connected compact riemannian manifold whose rational cohomology algebra 
 requires at least two generators then its free loop space homology has unbounded homology 
 and hence (Gromoll-Meyer \cite{MR0264551}) $M$ has infinitely many geometrically distinct closed geodesics. 

\medskip
  For completeness we make a brief review of the state of the art of the Sullivan 
  approach extension  to non-nilpotent spaces. 
  In fact,  Sullivan  minimal model is defined for any connected space 
  and it contains a lot of interesting properties  (\cite{FHT2}).
Moreover, a Sullivan model for a space $X$ with finitely many components $X_i$, $i = 1, \dots , n$, 
can be deduced easily from the models of the components. 
Denote by  $(M_Y,d)$   the minimal Sullivan model of a connected space $Y$. 
Then for $i = 1, \dots, n$, let   $a_i\in A_{PL}^0(X)$ be the $0$-form that takes the value $1$ on each point in $X_i$ 
and $0$ elsewhere. The $a_i$ form a family of orthogonal idempotents,  $1= \sum a_i$ 
and $a_i\cdot a_j = \delta_{ij}$ for all $i,j$, and  the direct sum 
$\displaystyle\bigoplus_i a_i\otimes (M_{X_i},d)$ is a Sullivan model for $X$.
Unfortunately the Sullivan model of a general non-nilpotent space 
does not give the rationalization of the higher homotopy groups.

\medskip
In the direction of models of non-nilpotent spaces, mention also that Bousfield and Kan define in \cite{bouskan}
other rational completions that $\Q_{\infty}X$. For instance (see Section VII.6.8 of op. cit.)  the completion
$\Q_{\infty}^*X$ has the same  fundamental group than $X$ and its universal covering is the $\Q$-completion of 
the universal covering of $X$. This has to be compared with the situation of $\Q_{\infty}X$ whose fundamental group
is the Malcev completion of the fundamental group of $X$.
In \cite{HalT}, the completion $\Q_{\infty}^*X$ has been algebraized
from Sullivan models of the fibrewise rationalization of a fibration. 
In the case of the universal covering  of a space $X$, these  models are ``rational vector bundles'' 
over $B\pi_{1}(X)$ with fiber a classical minimal model of the universal covering space, 
together with a homotopy class of a connection.
One asks if the same process gives a Lie version of this partial completion. The main advantage of such
Lie model should be to remove the hypothesis of finite type, required as usual in Sullivan's theory.

\medskip
Our program is carried out in the next sections, whose headings are self-explanatory.
\tableofcontents

\section{Complete differential graded Lie algebras}\label{sec:cdgl}

 Any algebraic object is considered $\bz$-graded and over  $\bq$.
 
 \begin{definition}\label{def:dglmc}
A \emph{differential graded Lie algebra} (dgl) is a graded vector space $L$ together with a 
\emph{Lie bracket}  $[\,\,,\,]$, satisfying graded antisymmetry and Jacobi identity, and a linear derivation   
$d$  of degree $-1$ such that $d^2= 0$. 

A \emph{Maurer-Cartan element} (or MC element) is an element $a\in L_{-1}$ such that $d a+ \frac{1}{2}[a,a] = 0$. 
We denote by ${\MC}(L)$ the set of Maurer-Cartan elements. Given $a\in\mc(L)$, 
the derivation $d_{a}=d+\ad_a$ is again a differential on~$L$. 

The  \emph{component} of $L$ at $a\in \mc(L)$ is the truncation of the perturbed $(L,d_a)$ at non-negative degrees,
$$L^{a}= L_{>0} \oplus \left( L_0\cap \ker d_a\right).
$$
 A dgl $L$ is called \emph{free} if it is free as a Lie algebra, that is, $L = \mathbb L(V)$ for some graded vector space $V$.
\end{definition}

\begin{definition} 
A \emph{topological dgl} is a dgl  with a topology making continuous the differential, the addition and the bracket. 
A \emph{morphism of topological dgl's} is a continuous map that is a dgl morphism. \end{definition}


We are interested in topological dgl's whose topologies come from filtrations.

\begin{definition}\label{def:dglfiltered}
A \emph{descending filtration} of a dgl  $L$ is a descending filtration   
$L = {\mathcal F}_1L \supset {\mathcal F}_2L \supset \cdots $ compatible with the differential and the bracket, $  d({\mathcal F}_kL )\subset {\mathcal F}_kL$ and $[{\mathcal F}_kL, {\mathcal F}_pL] \subset {\mathcal F}_{k+p}L$.
The filtered dgl $(L,{\mathcal F}_kL)$ is called \emph{complete} if
$$L = \varprojlim_k L/{\mathcal F}_kL\,.$$
In this case, we say also that ${\mathcal F}_kL$ is a \emph{complete descending filtration} on $L$.

 The \emph{completion} $(\widehat{L},{\mathcal F}_{k}\widehat{L})$ of a filtered dgl $(L,{\mathcal F}_kL)$ is given by
$$\widehat{L} = \varprojlim_n L/L^n\,\text{ and }
{\mathcal F}_{k}\widehat{L}=\left\{ \sum_{r\geq k}^{\infty}a_{r}\mid 
a_{r} \in {\mathcal F}_{r}L\right\}.$$
\end{definition}

For any dgl, $L$, we can consider the lower central series $L^1=L$, 
and for $n\geq 2$, $L^n = [L, L^{n-1}]$. If there is no reference to a filtration, the completion 
$\widehat{L}=\varprojlim_{n}L/ L^n$ of $L$
means the completion of $L$ for the lower central series.
The completion $\widehat{L}$ is then a topological dgl in which the open neighbourhoods of 0 are the subspaces 
$\ker(\widehat{L}\to L/L^n)$.


An element $\overline{a}$ of $\widehat{L}$ is a sequence $\overline{a}= (a_1, a_2, \dots)$ 
with $a_i\in L/L^i$ and $a_i =a_{i-1}$ in $L/L^{i-1}$.  
In particular, $L$ is dense in $\widehat{L}$, and $\widehat{L}$ is a complete dgl.

\begin{example}  
We write $\widehat{\mathbb L} (V)   = \widehat{\mathbb L(V)}$. Each element of $\widehat{\mathbb L}(V)$ can be seen as a series $\sum_n x_n$ with $x_n \in {\mathbb L}^n(V)$ for all $n$. This element is the limit of the sequence 
$\left(\sum_{n=1}^p x_n\right)_{p}$. 
Therefore if $f\colon (\widehat{\mathbb L}(V),d)\to (L,d)$ is a morphism of complete continuous dgl's, 
then $f(\sum_{n} x_n) = \sum_{n} f(x_n) $.
\end{example}

\begin{definition} We denote by {\bf cDGL} the category of   complete topological dgl's. Limits are the same than in the category of dgl's. The colimit of a diagram is the   completion of its colimit in the category of dgl's. In particular, the coproduct in $\catcdgl$ is denoted by $\widehat\amalg$. For sake of simplicity in the future we will say cdgl for topological complete dgl and morphism for continuous morphism.\end{definition}


\begin{definition}\label{def:BCH}
Given a cdgl $L$, the Baker-Campbell-Hausdorff product $*$ equips the vector space $L_0$,
of elements of degree~0,  with a group structure. The {\em gauge action} of $(L_0,*)$ on $\mc(L)$ is defined by
$$
x\mathcal G a=e^{\ad_x}(a)-\frac{e^{\ad_x}-1}{\ad_x}(dx).
$$
We denote by $\widetilde\mc(L)=\mc(L)/\mathcal G$ the set of equivalence classes of Maurer-Cartan elements modulo the gauge action.
\end{definition}

When $L$ is a cdgl and $A$ is a cdga, the complete tensor product
$L\widehat{\otimes} A=\varprojlim_{n}L/L^n\otimes A$
is a cdgl  whose bracket is defined by
$$[l\otimes a,l'\otimes a']=(-1)^{|a|.| l'|}[l,l']\otimes aa'.$$
The Deligne-Getzler $\infty$-groupoid associated to $L$ is the simplicial set 
$$\MC_{\bullet}(L)=\MC(L\widehat{\otimes} A_{PL}(\Delta^{\bullet}))$$
where $A_{PL}(\Delta^n)$ is the cdga of PL-forms on the simplex $\Delta^n$. As shown by M. Schlessinger and J. Stasheff
(\cite{SSth}), we have $\pi_{0}(\MC_{\bullet}(L))=\widetilde\mc(L)$.

In \cite{DR}, Dolgushev and Rogers establish a connection between quasi-isomorphisms and Maurer-Cartan elements.
We recall this result.

\begin{theorem} \label{propDR} \cite[Theorem 1.1]{DR}. Let $L$ and $L'$ be  cdgl's equipped with complete descending filtrations. If a morphism $f \colon L\to L'$  preserves the filtrations and
 induces  for all $m$ a quasi-isomorphism ${\mathcal F}_mL\to {\mathcal F}_mL'$, then, 
the simplicial map $\MC_{\bullet}(f)$ is a weak equivalence. In particular, the map
$\widetilde\mc(f)\colon \widetilde\mc(L)\xrightarrow{\cong}\widetilde\mc(L')$ is a bijection.
\end{theorem}

\section{The model and realization functors}\label{sec:modreal}

For each $n\ge 0$, let $(C_*(\Delta^n),\delta)$ be the usual simplicial chain complex on the standard $n$-simplex. 
We denote by $a_{i_{0}\dots i_{k}}$ the generator of the complete graded Lie algebra
$\libc(s^{-1}C_{*}(\Delta^n))$
associated to the $k$-face $\langle i_{0},\dots,i_{k}\rangle$.

The coface map $\sigma_{j}\colon \libc(s^{-1}C_{*}(\Delta^n)) \to \libc(s^{-1}C_{*}(\Delta^{n-1}))$
is then defined by
$\sigma_{j}(a_{i_{0}\dots i_{k}})=a_{r_{0}\dots r_{k}}$,
with $r_{s}=i_{s}$ for $s\leq j$ and $r_{s}=i_{s-1}$ for $s>j$,
and the convention $a_{r_{0}\dots r_{k}}=0$ if two indices coincide.

 In \cite{pri} we proved:

\begin{theorem}\label{bueno}{\em \cite[Theorems 2.3, 2.8 and 3.4]{pri}}
 there is a unique (up to isomorphism) cdgl of the form ,
$$\lasu_{\Delta^n}=\bigl(\libc\bigl(s^{-1}C_*(\Delta^n)\bigr),d\bigr),
$$ such that:
 \smallskip

 (1) For each $i=0,\dots,n$, the generators $a_i\in s^{-1}C_0(\Delta^n)$, corresponding to vertices, are Maurer-Cartan elements.
\smallskip

 (2) The linear part $d_1$ of $d$ is precisely  the desuspension of $\delta$.
\smallskip

(3)  The maps induced by the cofaces are the previous  $\sigma_{j}\colon \lasu_{\Delta^{n}}\to\lasu_{\Delta^{n-1}}$\end{theorem}

We briefly recall this construction. 
For a subsimplicial complex $X\subset \Delta^n$, we denote by $\lasu_{X}$ the sub-cdgl of $\lasu_{\Delta^n}$
generated by the elements $a_{i_{0}\dots i_{k}}$ corresponding to the simplices in $X$. As usual,
$\partial \Delta^n$ is the boundary of $\Delta^n$ and for each $i=0,\dots,n$,
$\wedge_{i}^n$
denotes the $i$-th horn of $\Delta^n$, obtained by removing the only $n$-simplex and the $(n-1)$-simplex 
$(0,\dots,\hat{\imath},\dots,n)$.

For $n=0$, $\lasu_{\Delta^0}$ is the free Lie algebra $\lib(a)$ generated by a Maurer-Cartan element.

For $n=1$,
$$ \lasu_{\Delta^1}=(\libc(a,b,x),d)\,,$$
is the  Lawrence-Sullivan interval.

 For $n=2$,
$$
\lasu_{\Delta^2}=(\libc (a_0,a_1,a_2,a_{01},a_{12},a_{02},a_{012}),d),
$$
in which $(\libc(a_0, a_1, a_{01}),d)$, $(\libc(a_1, a_2, a_{12}),d)$ and $(\libc(a_0, a_2, a_{02}),d)$ are  Lawrence-Sullivan intervals and
$$
d_{a_0}a_{012}=a_{01}*a_{12}*a_{02}^{-1}.
$$

For $n\ge 3$, $\lasu_{\Delta^n}$ is constructed inductively so that
$$
d_{a_0}a_{0\dots n}\in\lasu_{\partial\Delta^n}
$$
and with two additional important features,
 $$
 H(\lasu_{\Delta^n},d_{a_0})=H(\lasu_{\wedge^n_i},d_{a_0})=0.
 $$
  Suppose $\lasu_{\Delta^{n-1}}$ has been built and let $x=d_{a_0}a_{0\dots n-1}\in\lasu_{\partial\Delta^{n-1}}\subset\lasu_{\wedge^n_n}$. Since $H(\lasu_{\wedge^n_n},d_{a_0})=0$, there exists $y\in\lasu_{\wedge^n_n}$ such that $x=d_{a_0}y$. We set,
   $$
   d_{a_0}a_{0\dots n}=\Omega\quad\text{where}\quad \Omega=(-1)^n(a_{0\dots n-1}-y).
   $$
   In particular,
  $$
  d_{a_0}a_{0\dots n}-(-1)^na_{0\dots n-1}\in \lasu_{\wedge^n_n},
  $$
 which readily implies the following.
 
 \begin{lemma}\label{landni} 
 For $n\geq 2$ and $i\leq n$, there is a cdgl isomorphism
$$
 \lasu_{\Delta^n}\stackrel{\cong}{\longrightarrow}\lasu_{\wedge^n_i}\,\widehat\amalg
 \,\lib(u,du)
 $$
 with $u$ in degree $n-1$.
 \end{lemma}
 The same decomposition occurs in the case $n=1$ with a specific argument.

\begin{theorem}\label{LST} 
The morphism
  $$\psi\colon \lasu_{\Delta^1}=(\libc(a,b,x),d)\to (\libc(a,u,su), d), \hspace{5mm} dsu = u,$$
  defined by
  $$\psi (a) = a\,, \psi (x) = su\,, \psi (b) =  e^{\ad_{-su}}(a)+\frac{e^{\ad_{-su}}-1}{\ad_{-su}}(u),$$
is an isomorphism of cdgl's.
\end{theorem}

\begin{proof}
Let $i$ be the degree +1 derivation on $\libc(a,u,su)$ given by
$i(a)=su$, $i(su)=i(u)=0$. For the degree 0 derivation
$\theta=i\circ d +d \circ i $, we get
  $\theta(u)=\theta(su)=0$, and
$$\theta(a)= (i\circ d +d \circ i)(a)=
d su-\frac{1}{2}i[a,a]=u-\ad_{su}(a)= u+ \ad_{-su}(a).$$
 We deduce by induction,
$$
\theta^n(a)
=
\ad_{-su}^{n-1}(u)+\ad_{-su}^n(a).
$$
Therefore
$$
  e^{\theta}(a)
=
  \sum_{i\geq 0}\frac{\ad_{-su}^{i}(a)}{i!}+
 \sum_{i\geq 0}\frac{\ad_{-su}^i(u)}{(i+1)!}
=
 e^{\ad_{-su}}(a)+\frac{e^{\ad_{-su}}-1}{\ad_{-su}}(u)= \psi (b).
$$
This formula can be reversed and it  gives
$$
 u
 =
\frac{\ad_{-su}}{e^{\ad_{-su}}-1} e^{\theta}(a)
-\frac{\ad_{-su}\,e^{\ad_{-su}}}{e^{\ad_{-su}}-1} (a).
$$
 In particular, we obtain the surjectivity of $\psi$ from
$$
u=
\psi\left(\frac{\ad_{-x}}{e^{\ad_{-x}}-1}(b-a)+\ad_{x}(a)\right)=\psi (d x),
$$
and $d\psi (x) = \psi (d x)$
 \end{proof}

\begin{remark}\label{rem:LSinterval}

(i)  The inverse isomorphism
$$\psi^{-1} \colon (\libc(a,u,su),d) \stackrel{\cong}{\longrightarrow} (\libc(a,b,x),d)$$
 is given by
 $  \psi^{-1} (a)=a$, $\psi^{-1} (su)= x$ and $\psi^{-1} (u) =  \partial x$.

 (ii) The differential of $\lasu_{\Delta^1}$ is entirely determined by the differentials of $a$ and $b$, and by
the linear part of the differential of $x$, see \cite[Theorem 1.4]{pri},  \cite[Main Theorem]{PPTD}.
\end{remark}

The construction of $\lasu_{X}$ can be extended to any simplicial set $X$. Let us identify, as usual,
any simplex $\sigma\in X_{n}$ with a simplicial map $\sigma\colon\underline{\Delta}^n\to X$, where
$\underline{\Delta}^n$ denotes the simplicial set whose $p$-simplices are integer sequences $0\le i_0\le\dots\le i_p\le n$. Then, $X$ can be recovered from its simplices as the colimit 
$$
X=\varinjlim_{\sigma\in X}\underline\Delta^{|\sigma|}.
$$

\begin{definition}
The {\em model} of any simplicial set $X$ is defined as the cdgl
$$
 \lasu_X=\varinjlim_{\sigma\in X}\lasu_{\Delta^{|\sigma|}}.
$$
\end{definition}

Whenever $X$ is a finite simplicial set, that is, a simplicial set with a finite number of non-degenerate simplices, its model is the free complete Lie algebra
$$
\lasu_X=\libc\bigl(s^{-1}NC_*(X)\bigr)
$$
generated by the non-degenerate simplices of $X$, and whose differential $d$ is completely determined by the following:

 (1) The non-degenerate $0$-simplices are Maurer-Cartan elements.
\smallskip

 (2) The linear part $d_1$ of $d$ is precisely  the desuspension of the differential in $NC_*(X)$.
\smallskip

(3) If $j\colon Y\subset X$ is a subsimplicial set, then
$
\lasu(j)=\libc\bigl(s^{-1}NC_*(j)\bigr)$.

\begin{example}
When $X$ is a finite simplicial complex, and $Y$ is the associated simplicial set,
then the above definition of $\lasu_{Y}$ and the definition of $\lasu_{X}$ as sub-cdgl of some $\lasu_{\Delta^n}$
coincide.
\end{example}

On the other hand  $\lasu_{\Delta^{\bullet}}=\{\lasu_{\Delta^n}\}_{n\ge 0}$ has a cosimplicial cdgl structure \cite[\S3]{pri} which gives rise to the following.

\begin{definition}\label{def:realization}
 The {\em realization} of a cdgl $L$ is defined as the simplicial set
$$
 \langle L\rangle=\Hom_{\catcdgl}(\lasu_{\Delta^\bullet},L).
$$
\end{definition}

In \cite[Theorem 4.1]{cuar} (cf. also \cite{ni1}), we prove that the realization $\langle L\rangle$  coincides 
with the Deligne-Getzler simplicial set $\MC_{\bullet}(L)$,  \cite{ber,getz,hi}.

\medskip
In \cite[Theorem 4.3]{pri}, we prove that the above constructions define a pair of adjoint functors
$$
\xymatrix{  \catss \ar@<0.75ex>[r]^-{\lasu} &\catcdgl. \ar@<0.75ex>[l]^(.54){\langle-\rangle} }
$$
For any cdgl $L$, its realization is the disjoint union \cite[Proposition 4.6]{pri},
\begin{equation}\label{componentes}
\langle L\rangle=
{\sqcup}_{a\in\widetilde\mc(L)}\bigl\langle L^{a}\bigr\rangle,
\end{equation}
and the homotopy groups of each of these components are \cite[Proposition. 4.5]{pri},
\begin{equation}\label{gruposhomotopia}
\pi_n \bigl\langle L^{a}\bigr\rangle\cong H_{n-1}\bigl(L^{a}\bigr),\quad n\ge 1,
\end{equation}
in which $H_0$ is considered with the group structure given by the Baker-Campbell-Hausdorff product.

On the other hand we proved in \cite{ter} that, given $X$ any finite simplicial set, $\widetilde\mc(\lasu_X)=\pi_0(X^+)$ where $X^+$ denotes the disjoint union of $X$ with a point. This, together with (\ref{componentes}), gives,
\begin{equation}\label{picero}
\pi_0\langle\lasu_X\rangle=\pi_0(X^+).
\end{equation}
Now, we determine the homotopy type of each of these components, and we begin with the component at 0, 
see Definition~\ref{def:dglmc}.

\begin{theorem}\label{componentecontractil} Given $X$ a finite 
simplicial set, then $H(\lasu_{X})=0$ and the realization of the component at $0$,
 $\bigl\langle\lasu_X^{0}\bigr\rangle$, is contractible.
\end{theorem}

\begin{proof}
 Write ${\lasu}_X = (\libc(V),d)$. We choose an element $a\in V_{-1}$ corresponding to a non-degenerate vertex.
 We denote by $W$ the subspace of $V$ defined by  $W_n = V_n$ for $n\geq 0$ and such that
 $W_{-1}$ is the subspace of $V_{-1}$  generated by the differences $b_i-a$
 where the $b_i$'s  browse the vertices of $X$.
 Then, we have an isomorphism of complete Lie algebras,
 $$\libc(V) \cong \L(a)  \,\widehat\amalg \, \libc(W).$$
 In ${\lasu}_X$, the differential $d$ satisfies
 $$da= -\frac{1}{2} [a,a] \quad\text{and}\quad dx = d_ax-[a,x]\,.$$
 Denote by $\mathcal I$ the ideal generated by $W$ in ${\lasu}_X$.
 Since $V$ is finite dimensional and $| a| = -1$, the ideal $\mathcal I$ is the complete free Lie algebra on the
 vector space $\overline{W}$, of basis $(\ad^n_a(v_i))_{i,n}$ where $(v_i)_{i}$ is a basis of $W$ and  $n\geq 0$, i.e.,
 $$\mathcal I = \libc(\overline{W}).$$
 We denote by  $d_1$ the linear part of the differential $d$ in $\lasu_{X}$ and
 by $\delta$ the linear part of the  differential induced by $d$ in  $\mathcal I=\libc(\overline{W})$.
 An induction shows that
 $$\delta (\ad_a^n (v)) =
 \renewcommand{\arraystretch}{1.3}
 \left\{\begin{array}{ll}
 (-1)^n \ad_a^n (d_1(v)) - \ad_a^{n+1}(v)\hspace{3mm}\mbox{} & \mbox{ if $n$ is odd,}\\
 (-1)^n \ad_a^n(d_1(v)) & \mbox{ if $n$ is even.}
 \end{array}\right.
 \renewcommand{\arraystretch}{1}$$
 From this computation, we deduce that if
 $\sum_{n,i} \alpha_{n,i}\, \ad_a^n (v_i)$, $\alpha_{n,i}\in \bq$,  is a $\delta$-cycle, then we have
 $\alpha_{2n+1,i}=0$, for any $n$ and $i$, and $\sum_{n,i}\alpha_{2n,i}\,\ad^{2n}(d_{1}v_{i})=0$. Thus, we obtain
$$\sum_{n,i}\alpha_{n,i}\, \ad_a^n (v_i) = -\delta( \sum_{n,i} \alpha_{2n,i}\, \ad_a^{2n-1}(v_i))\,.$$
Therefore,  $H(\mathcal I, \delta) = 0$ and from
 \cite[Proposition 2.4]{pri}, we deduce that $\cI$ is $d$-acyclic.
Finally, since ${\lasu}_X/\mathcal I = (\lib(a),d)$, we get $H({\lasu}_X)= 0$.

Since $\lasu_{X}$ is acyclic, $H(\lasu^0_{X})=0$.
By equation (\ref{gruposhomotopia}) this is equivalent to  $\pi_*\bigl\langle\lasu_X^{0}\bigr\rangle=0$ and the theorem follows.
\end{proof}

Consider now a connected, finite simplicial  set, $X$ and a 0-simplex, $a$, of $X$. The cdgl $(\lasu^a_{X},d_{a})$
is then directly related to the Sullivan minimal model of $X$, $(\Lambda V,d)$.

To specify this relation, let us recall adjoint functors, introduced by D. Sinha and B. Walter
in \cite{SinWal}, and that are a dual version of the Quillen adjunction $(\mathcal L, C^*)$.
First, let us remind that a graded Lie coalgebra is a graded vector space $E$
together with a comultiplication 
$\Delta\colon E\to E\otimes E$ that satisfies two properties
$$(1+\tau)\circ \Delta=0\quad \text{and}\quad 
(1+\sigma+\sigma^2)\circ (1\otimes \Delta)\circ \Delta=0,$$
where $\tau\colon E\otimes E\to E\otimes E$ is the graded transposition and
$\sigma\colon E\otimes E\otimes E \to E\otimes E\otimes E$ the graded cyclic permutation.
For instance, the tensor coalgebra on a graded vector space $V$, $T^c(V)$, admits a graded Lie coalgebra structure 
defined by $\Delta_{L}=\Delta-\tau\circ\Delta$. The free lie coalgebra on $V$, $\L^c(V)$,
is the quotient of $T^c(V)$ by the set of decomposable elements for the shuffle product.
Let us denote by \\
-- $\catcdga$ the category  of  commutative, augmented, differential $\Z$-graded algebras (cdga for short),\\
-- $\bf DGLC$ the category  of differential $\Z$-graded Lie coalgebras (dglc for short).

\medskip
In \cite{SinWal}, D. Sinha and B. Walter introduce adjoint functors
$$
\xymatrix{  {\bf DGLC} \ar@<0.75ex>[r]^-{\mathcal A} &\catcdga \ar@<0.75ex>[l]^(.54){\mathcal E} }
$$
and establish their properties with some restriction of connectivity. We use them in the $\Z$-graded case and
give in Lemma~\ref{lem:SinhaWalter} a proof of the  properties we need in this context.
First, these functors are defined as follows.\\
-- If $(E,d)$ is a dglc of comultiplication $\Delta$, we set
$${\mathscr A}(E,d)=(\Lambda s^{-1}E,D)$$
where $\Lambda s^{-1}E$ is the free commutative graded algebra on the suspension of $E$ and
$$Ds^{-1}x=\frac{1}{2}\sum_{i}(-1)^{|x_{i}|} s^{-1}x_{i}\Lambda s^{-1}x'_{i}-s^{-1}dx,$$
if $\Delta x=\sum_{i} x_{i}\otimes x'_{i}$.
In particular, if $E$ is a finite dimensional cdgl, then ${\mathcal A}(E)= C^*(E^\#)$, where $E^\#$ denotes the dual dgl
of $E$. \\
-- Let $A$ be a cdga with ideal of augmentation $\ov{A}$. Observe that the differential of the Bar construction,
$BA=T^c(s\ov{A})$,
preserves the ideal generated by the shuffle products. Therefore, it induces a differential $d$ on
$\L^c(s\ov{A})$ and we set
$${\mathscr E}(A)=(\L^c(s\ov{A}),d).$$

 \begin{lemma}\label{lem:SinhaWalter}\mbox{}
 \begin{enumerate}
\item  The functor ${\mathcal E}$ preserves the quasi-isomorphisms.
\item If $f\colon E\to E'$ is a dglc quasi-isomorphism such that $E=E_{\geq 0}$
and $E'=E'_{\geq 0}$, then 
${\mathscr A}(f)$ is a  quasi-isomorphism.
\item If  $A$ is a cdga such that $\ov{A}=\ov{A}^{\geq 1}$, then the adjunction map 
 ${\mathcal A}{\mathcal E}(A)\to A$ 
 is a quasi-isomorphism.
 \end{enumerate}
  \end{lemma}

\begin{proof}
(1) For any  cdga $A$, the dglc ${\mathcal E}(A)= (\mathbb  L^c(s\overline{A})$,d) 
is the union of the sub dglc's 
$(({\mathbb L}^c)^{\leq n}(s\overline{A}),d)$, 
image of $T^{\leq n}(s\overline{A})$ along the projection $BA\to {\mathcal E}A$. 
Let $f \colon A\to A'$ be a cdga quasi-isomorphism. 
By induction on $n$ and the five lemma, the map $(\mathbb L^c)^{\leq n}(f)$ is a quasi-isomorphism. 
This implies the result.

(2) We filter ${\mathscr A}(E)=(\Lambda sE,D)$ and ${\mathscr A}(E')=(\Lambda sE',D)$ by the length of words in $\Lambda$.
At  level~0 of the spectral sequences, the map ${\mathscr A}(f)$ induces  
$\Lambda sf\colon \Lambda (sE, -sd)\to \Lambda (sE',-sd)$ and $E^1(\cA(f))=\Lambda sH_{*}(f)$ which is an isomorphism.
By hypothesis, all elements are in degrees $\geq 1$, thus the spectral sequences converge and ${\mathscr A}(f)$
is a quasi-isomorphism.
 
 (3) $\bullet$ \emph{We first suppose that $A$ is finitely generated.}
  Then   ${\mathcal E}(A)$ is infinite countable. 
  The elements of  $\Lambda^p(s^{-1}{\mathcal E}(A))$ are said of length $p$
  and we define   $F_p{\mathcal A}{\mathcal E}(A)$ as the subcomplex generated by the elements $\alpha$
   with degree ($\alpha$) - length ($\alpha$) $\geq p$. 
   We introduce also a decreasing filtration on $A$, by $F^0A= A$ and $F^1A= 0$. 
   The adjunction map ${\mathcal A}{\mathcal E}(A)\to A$ is a morphism of filtered complexes, and since $\overline{A}= \overline{A}^{\geq 1}$, the induced spectral sequences converge.
 
 The term $(E^0({\mathcal A}{\mathcal E}(A)), d_0)$ is isomorphic to the functor ${\mathcal A}$ applied 
 to the free Lie coalgebra $\mathbb L^c(s\overline{A})$ with  differential $0$. 
  Denote  by $E $ a finite type graded vector space satisfying $ {E} = {E}^{\geq 1}$ and isomorphic to $s\overline{A}$ 
  as a (non graded) vector space. We then have an isomorphism of (nongraded) complexes 
  ${\mathcal A}(\mathbb L^c(E), 0) \cong C^*(\mathbb L(E^\#), 0= C^*{\mathcal L}^*(\mathbb Q \oplus E)$, 
  where $\mathbb Q\oplus E$ is equipped with the trivial multiplication. 
  Since the map $C^*{\mathcal L}^*(\mathbb Q\oplus E)\to (\mathbb Q \oplus E, 0)$ is a quasi-isomorphism, 
  the original map $(E^0({\mathcal A}{\mathcal E}(A), d_0) \to (A,0)$ is a quasi-isomorphism. The convergence of the spectral sequences implies then the result.
 
 $\bullet$ \emph{We consider now the general case.}
  Each cdga $A$ is the union of finitely generated sub cdga's $A_\alpha$, 
  and the adjunction map ${\mathcal A}{\mathcal E}(A)\to A$ is the limit of the adjunction maps 
  $ {\mathcal A}{\mathcal E}(A_\alpha)\to A_\alpha$. The result comes from the commutativity of the homology
  with inductive limits.
 \end{proof}
 
 \begin{proposition}\label{prop:twomodels}
 Let $X$ be a connected finite simplicial set. Denote by $(\Lambda V,d)$
 its minimal Sullivan model and by
$L= (\lasu_X^a,d_a)$ the component of the Lie model $\lasu_{X}$ at the Maurer-Cartan element $a$.
Then there exist quasi-isomorphisms,
$$(\Lambda V,d)
\simeq
\varinjlim_{n} C^*(L/L^n).
$$
 \end{proposition}
 
 \begin{proof}
 From  \cite[Theorem 7.4]{pri} 
 we know the existence of a  dglc $L^c_X $, quasi-isomorphic to $\mathcal E(A_{PL}(X))$ and whose dual 
is the cdgl  $ (\lasu_X^a,d_a)$.
From \cite[Corollary 5.4]{pri}, we may replace $ (\lasu_X^a,d_a)$ by a cdgl $L=(\hL(sH_{*}(X;\Q),d)$
and therefore, as $sH_{*}(X;\Q)$ is a finite type vector space, we may replace $L^c_X $
by $M^c_{X}=(\L^c(sH^*(X;\Q),d)$.
The dglc $M^c_X$ can be written as  the union of finite dimensional dglc's $M_X^c(n)$,
 whose duals are the quotients $L/L^n$. 

From the compatibility of $\mathscr E$ with quasi-isomorphisms, we deduce
$${\mathscr E}(\Lambda V,d)\simeq {\mathscr E}(A_{PL}(X))\simeq L^c_{X}\simeq M^c_{X}.$$

From the property of the adjunction map, we have
$$(\Lambda V,d)\simeq {\mathscr A}{\mathscr E}(\Lambda V,d).$$

Finally, the compatibility of $\mathscr A$ with quasi-isomorphisms can be used here and we get
$${\mathscr A}{\mathscr E}(\Lambda V,d)\simeq {\mathscr A}(M^c_{X}).$$
Since
${\mathcal A}$ is a left adjoint functor, it is compatible with inductive limits and we  obtain the desired quasi-isomorphism,
$$
(\Lambda V,d) \simeq {\mathcal A}(\varinjlim_n M_X^c(n)) \simeq \varinjlim_n {\mathcal A} (M_X^c(n)) = \varinjlim_{n} C^*(L/L^n)\,.
$$
 \end{proof}

For the non trivial components the realization of $\lasu_X^a$ is related to the Sullivan realization 
$\langle \Lambda V,d\rangle_S$ of the Sullivan minimal model and to the Bousfield-Kan completion of $X$.

\medskip
We recall that in \cite{bouskan}, A. Bousfield and D. Kan define a $\bq$-completion functor  that associates to each
simplicial set  $X$, a simplicial set $\bq_{\infty}X$ together with a natural morphism,
$X\to\bq_{\infty}X$. This construction is governed by the following property:
a map $f\colon X\to Y$ induces an isomorphism in rational homology,
$\widetilde H_{*}(X;\bq)\stackrel{\cong}{\to}\widetilde H_{*}(Y;\bq)$,
if and only if
the map
$\bq_{\infty}f\colon \bq_{\infty}X\to \bq_{\infty}Y$
is a homotopy equivalence. The space $\bq_{\infty}X$ is called the \emph{$\bq$-completion of $X$.}
Moreover, if $(\Lambda V,d)$ denotes the Sullivan minimal model of $X$, then A. Bousfield and 
V.  Gugenheim 
(\cite[Theorem 12.2]{bousgu}) prove that the Sullivan realization  $\langle \Lambda V,d\rangle_S$ is 
homotopy equivalent to the Bousfield-Kan completion of $X$. 

\medskip
On the other hand, the
\emph{rationalization} of a nilpotent simplicial set $X$ is a simplicial set $X_{\bq}$ equipped with a morphism
$\rho\colon X\to X_{\bq}$ satisfying the following:

(i) The reduced homology $\widetilde H_{*}(X_{\bq};\bz)$ is a rational vector space.

(ii)  The induced map $H_{*}(\rho;\Q)$ is an isomorphism.\\
These two properties characterize the rationalization up to a homotopy equivalence. 
For a nilpotent simplicial set, $X$, the $\Q$-completion and the rationalization coincide,
$\bq_{\infty}X\simeq X_{\Q}$, see
\cite[Chapter V, 4.3]{bouskan}.

\begin{theorem}\label{componentecompleccion} 
Let $X$ be a connected finite simplicial set  and let $a\in \lasu_X$ be a non-trivial Maurer-Cartan element. 
Then, the realization of $\lasu_X^{a}$ is homotopy equivalent to the $\Q$-completion of $X$,  
$\bigl\langle \lasu_X^{a}\bigr\rangle\simeq \bq_\infty X$.
\end{theorem}

\begin{proof} 
Let $L = \lasu_X^a$ and $(\Lambda V,d)$ the Sullivan minimal model of $X$. By Proposition~\ref{prop:twomodels}, 
$(\Lambda V,d)\simeq \varinjlim_n C^*(L/L^n)$. Since the Sullivan realization functor is adjoint to the PL-forms functor, we 
deduce:
$$\langle \Lambda V,d\rangle_S \simeq \varprojlim_n \langle C^*(L/L^n)\rangle_S\,.$$
As $L/L^n$ is finite dimensional,  the Getzler isomorphism (\cite[Proposition 1.1]{getz}) gives
a weak-equivalence of simplicial sets,
$$\langle C^*(L/L^n)\rangle_S \simeq \MC_{\bullet}(L/L^n)\,.$$
Now by \cite{ber}, (see also \cite{cuar} or \cite{ni1}),
the simplicial set $\MC_{\bullet}(L/L^n)$ is weakly equivalent to $\langle L/L^n\rangle $. 
Finally, since the realization functor of  cdgl's is a right adjoint, it commutes with direct limits and
we have:
$$\langle \Lambda V,d\rangle_S\simeq \varprojlim_n \langle C^*(L/L^n)\rangle_S 
\simeq \varprojlim_n \langle L/L^n\rangle = \langle\varprojlim_n L/L^n\rangle = \langle L\rangle \,.$$
The result follows since the Sullivan realization model is already weakly equivalent to the Bousfield-Kan completion.
\end{proof}

From Theorem~\ref{componentecompleccion} and equation (\ref{gruposhomotopia})
before Theorem~\ref{componentecontractil}, we deduce:

\begin{corollary}\label{cor:liehomotopy}
Let $X$ be a connected finite simplicial complex. The homotopy Lie algebra  $\pi_{*}(\bq_{\infty}X)\cong H_{*-1}(\lasu_{X}^{a})$ is isomorphic to
the homotopy Lie algebra of the Sullivan model of $X$.
\end{corollary}

\begin{remark}\label{prop:sullivannotfinite}
Let $X$ be a simply connected simplicial set whose Betti numbers are countable
and at least one of them is infinite, and let $(\Lambda V,d)$ be its Sullivan minimal model. Then,
$\langle \Lambda V,d \rangle_{S}\neq \bq_{\infty}X$.

We first observe that $\bq_{\infty}X=X_{\Q}$. For that, we write $X=\cup_{n}X_{n}$ as an increasing union of simply connected finite simplicial sets. The injection
$\iota\colon X\to \cup_{n}(X_{n})_{\Q}$
induces an isomorphism in rational homotopy, therefore $\Q_{\infty}(\iota)$
is a homotopy equivalence. The space $\cup_{n}(X_{n})_{\Q}$ being a simply connected space, it is  $\Q$-complete. 

Now let $p$ be the smallest degree $i$ such that $H^i(X;\bq)$ is infinite dimensional and decompose
$(\Lambda V^{\leq p},d)=(\Lambda V^{<p}\otimes \Lambda V^p,d)$.
As $\dim H(\Lambda V^{<p},d)<\infty$ and $\dim H^p(\Lambda V,d)=\infty$, we have $\dim V^p=\infty$.
Hence, $\dim \pi_{p} \langle \Lambda V,d\rangle_S =\dim (\Hom(V^p,\bq))$ is not countable.
Therefore,  $\dim H_{p} \langle \Lambda V,d\rangle_S $ is not countable and the assertion follows.
\end{remark}

\begin{example}[\emph{The circle}]
By viewing the circle as the boundary of the triangle, we get
$$\lasu_{S^1}=(\hL (a_{0},a_{1},a_{2},a_{01},a_{12},a_{02}),d),$$
with the notation and the differential of Lawrence and Sullivan intervals.
We use this example to show how non-trivial homotopy types or, equivalently, non-trivial homology of models,
arise when perturbing the differential by a Maurer-Cartan element.

First, for the differential $d_{a_{0}}$, we know that 
$$H_{0}(\lasu_{S^1},d_{a_{0}})=\Q(a_{01}\ast a_{12}\ast (a_{02})^{-1}),$$
and the map $(\L(u),0)\to (\lasu_{S^1},d_{a_{0}})$, sending the element $u$ of degree 0 on
$a_{01}\ast a_{12}\ast (a_{02})^{-1}$
is a quasi-isomorphism. 

On the other hand, for the differential $d$, we have from Theorem~\ref{componentecontractil},
$H_{0}(\lasu_{S^1},d)=0$. We can also observe that the element $a_{01}\ast a_{12}\ast (a_{02})^{-1}$ having 
the bracket $-[a_{0},a_{01}\ast a_{12}\ast (a_{02})^{-1}]$ as boundary is not a $d$-cycle.
\end{example}

\section{The model category structure on {\bf cDGL}}\label{sec:closedmodel}

Henceforth, by  {\em model category} we mean the original  closed model category definition of Quillen \cite{qui2}. We prove:

\begin{theorem}\label{model2} There is a cofibrantly generated model category structure on $\catcdgl$ for which:
\begin{enumerate}
\item[$\bullet$]
A morphism $f\colon A\to B$ is a fibration if 
it is surjective in non-negative degrees.
\item[$\bullet$]
A morphism $f\colon A\to B$ is a weak equivalence if $\widetilde\mc(f)\colon \widetilde\mc(A)\stackrel{\cong}{\to}\widetilde\mc(B)$ is a bijection and  $f^{a}\colon A^{a}\stackrel{\simeq}{\to} B^{f(a)}$ is a quasi-isomorphism for every $a\in\widetilde\mc(A)$.
\item[$\bullet$]
 A morphism is a cofibration if it has the left lifting property with respect to trivial fibrations.
 \item[$\bullet$] The cdgl morphisms
$$
\lasu_{\mathscr I}=\{\lasu_{\partial\Delta^n}\hookrightarrow \lasu_{\Delta^n}\}_{n\ge0},\quad
\lasu_{\mathscr J}=\{\lasu_{\wedge^n_i}{\hookrightarrow} \lasu_{\Delta^n}\}_{n\ge0, \,i=0,\dots,n},
$$
are  generating sets of cofibrations and trivial cofibrations respectively.
\end{enumerate}
 \end{theorem}

For the proof, we use, like in  the work of Bandiera \cite[Theorem 5.2.28]{ban}   the  version of the so called 
\emph{Transfer Principle}  in \cite[Sections 2.5 and  2.6]{bermoer} which is encoded in the next definition 
and the following theorem.

\begin{definition}  Let $\mathcal C$ be a  model category cofibrantly generated by the sets $\mathscr I$  and $\mathscr J$ of generating cofibrations and generating trivial cofibrations respectively. Let $\mathcal D$ be a category  with finite limits and small colimits, and let
$$
\xymatrix{  \mathcal C \ar@<0.75ex>[r]^-{F} &\mathcal D \ar@<0.75ex>[l]^(.52){G} }
$$
be a pair of adjoint functors (upper arrow denotes left adjoint).  Define a map
$f$ in $\mathcal D$ to be  a weak equivalence (resp. fibration) if $G(f)$ is a weak equivalence (resp. fibration). \end{definition}

\begin{theorem}\label{principle} The transfer principle defines a model category in $\mathcal D$ cofibrantly generated by the families $F(\mathscr I)$ and $F(\mathscr J)$ provided:

(i) The sets $F(\mathscr I)$ and $F(\mathscr J)$ permit the small object argument.

(ii) $\mathcal D$ has a functorial fibrant replacement and a functorial path object for fibrant objects. \end{theorem}

\medskip
We apply the Transfer Principle to the adjunction
$$
\xymatrix{  \catss \ar@<0.75ex>[r]^-{\lasu} &\catcdgl. \ar@<0.75ex>[l]^(.54){\langle-\rangle} }
$$
 We endow $\catss$ with the usual structure of closed model category:  
weak equivalences are the maps whose geometric realization is a weak equivalence,
 fibrations are the Kan fibrations
 and cofibrations are the injective maps.

\medskip
Thus, a   morphism $f$ of cdgl's is by definition a fibration (resp. a weak equivalence) 
if $\langle f \rangle$ is a fibration (resp. a weak equivalence) of simplicial sets. 
Weak equivalences and fibrations admit nice characterizations.

\begin{proposition}\label{we} The realization $\langle f\rangle$ of a cdgl morphism $f\colon A\to B$ is a weak equivalence if and only if $\widetilde\mc(f)$ is a bijection and $f^a\colon A^a\stackrel{\simeq}{\to}B^{f(a)}$ is a quasi-isomorphism for each $a\in\widetilde\mc(A)$.
\end{proposition}

\begin{proof} It follows immediately from the equalities (\ref{componentes}), (\ref{gruposhomotopia})
after Definition~\ref{def:realization}.
\end{proof}

\begin{proposition}\label{fib} The realization $\langle f\rangle$ of a cdgl morphism $f\colon A\to B$ is a fibration if and only if 
the restriction $f\colon A_{\geq 0}\to B_{\geq 0}$ is surjective.
\end{proposition}

\begin{proof}
  It   follows from Lemma \ref{landni} and Theorem \ref{LST} that for any cdgl $L$, any $n> 0$ 
  and any $i=0,\dots n$, the map $ f\mapsto \bigl(f_{|_{\lasu_{\wedge^n_i}}},f(a_{0\dots n})\bigr)$,
\begin{eqnarray}\label{iiso}
\mbox{\rm Hom}_{\catcdgl}(\lasu_{\Delta^n},L)\stackrel{\cong}{\longrightarrow}
\mbox{\rm Hom}_{\catcdgl}(\lasu_{\wedge^n_i},L)\times L_n,
\end{eqnarray}
is a bijection. In particular, given any $x\in L_n$, any cdgl morphism $f\colon \lasu_{\wedge^n_i}\to L$ extends to a cdgl morphism $f\colon\lasu_{\Delta^n}\to L$ such that $f(a_{0\dots n})=x$.

By definition, the realization
$\langle f \rangle\colon \langle A\rangle \to \langle B\rangle$
is a fibration if and only if there exists a dotted lifting in any square commutative diagram of the form,
$$\xymatrix{
 \wedge^n_{i}\ar[r]\ar[d]
 &
 \langle A\rangle \ar[d]^-{\langle f\rangle}\\
 \Delta^n\ar[r]
 \ar@{-->}[ru]
 &
  \langle B\rangle.
  }$$
  By adjunction, this is equivalent to the existence of the dotted lifting  in the corresponding diagram,
  $$\xymatrix{
  \lasu_{\wedge^n_{i}}\ar[r]\ar[d]&
  A\ar[d]^f\\
  \lasu_{\Delta^n}\ar[r]
  \ar@{-->}[ru]
  &
  B.
  }
  $$
  By (\ref{iiso})  the morphism $\lasu_{\Delta^n}\to B$ is uniquely determined  by its image 
  $x\in B_n$ of $a_{0\dots n}$, as it is $\lasu_{\wedge^n_i}\to A\stackrel{f}{\to} B$ on $\lasu_{\wedge^n_i}$. 
  Hence, if $f$ is surjective in non-negative degrees, choose $f(y)=x$ and define the lifting
  $\lasu_{\Delta^n}\to A$
  as the only morphism extending 
  $\lasu_{\wedge^n_i}\to A$
  and sending $a_{0\dots n}$ to $y$.

  Conversely, assume the dotted lifting exists for any such commutative square  and let 
  $x\in B_n$, $n\geq 0$. By (\ref{iiso})  there exists a unique cdgl morphism $\lasu_{\Delta^n}\to B$ which is zero on 
  $\lasu_{\wedge^n_i}$ and send $a_{0\dots n}$ to $x$. Hence, its lifting $\lasu_{\Delta^n}\to A$ sends $a_{0\dots n}$ to $y$ with $f(y)=x$.
\end{proof}

\medskip \begin{proof}[Proof of Theorem \ref{model2}] It remains to see that the model functor and the category $\catcdgl$ satisfy the conditions (i) and (ii) of Theorem \ref{principle}. First of all, observe that
$\catcdgl$ has   arbitrary limits and colimits. Note also that every cdgl is fibrant so the first assertion of (ii) is trivially satisfied. For the second, recall that a path object for $A$ is a factorization of its diagonal
$A\stackrel{\simeq}{\to} A^I \twoheadrightarrow A\times A$ into a
weak equivalence followed by a fibration. In $\catcdgl$ a functorial path object is given by   the sequence,
$$
\xymatrix{ L\ar@{^{(}->}[r]^(.27)\simeq & L^I=L \widehat
\otimes \Lambda(t,dt)
\ar@{->>}[rr]^-{(\varepsilon_0,\varepsilon_1)}&&
L\times L}
$$
with $|t|=0$, $|dt|=-1$. The fact that $L\hookrightarrow L^I$
is a weak equivalence is a direct consequence of Theorem~\ref{propDR} for the filtration
$\cF_{m}L^I=\cF_{m}L\otimes \Lambda(t,dt)$.

 For (i) recall that $\catss$ if cofibrantly generated by the generating sets,
$$
\mathscr I=\{\partial\Delta^n\hookrightarrow \Delta^n\}_{n\ge0},\quad
\mathscr J=\{\wedge^n_i\stackrel{\simeq}{\hookrightarrow} \Delta^n\}_{n\ge0, \,i=0,\dots,n}
$$
of cofibrations and trivial cofibrations respectively. Since a countable union of countable ordinals remains countable,  the  models
$
\{\lasu_{\partial\Delta^n}\}_{n\ge0}$  and\break $\{\lasu_{\wedge^n_i}\}_{n\ge0, \,i=0,\dots,n}$
satisfy the small object argument with respect to the morphisms
$$
\lasu_{\mathscr I}=\{\lasu_{\partial\Delta^n}\hookrightarrow \lasu_{\Delta^n}\}_{n\ge0},\quad
\lasu_{\mathscr J}=\{\lasu_{\wedge^n_i}{\hookrightarrow} \lasu_{\Delta^n}\}_{n\ge0, \,i=0,\dots,n}.
$$
\end{proof}

Note that the path object is brought in from the classical setting both in the commutative
\cite[\S 11]{FHTI} 
and the Lie 
\cite[Section II.5]{tan} sides.

\medskip
The main advantage of having a model structure obtained by the transfer principle is its compatibility with
the homotopy relation.

\begin{corollary} The realization and model functors,
$$
\xymatrix{  \catss \ar@<0.75ex>[r]^-{\lasu} &\catcdgl, \ar@<0.75ex>[l]^(.52){\langle-\rangle} }
$$
form a Quillen pair. In particular, they preserve weak equivalences and induce adjoint functors in the homotopy categories,
$$
\xymatrix{  \ho\catss \ar@<0.75ex>[r]^-{\lasu} &\ho\catcdgl. \ar@<0.75ex>[l]^(.54){\langle-\rangle} }
$$
\end{corollary}

\begin{corollary}
\mbox{}
\begin{enumerate}
\item[(i)] For any cylinder in $\catcdgl$, the realization of any two homotopic maps with cofibrant domain,
 $f,g\colon A\to B$,  are homotopic maps of simplicial sets.
\item[(ii)] The model of two homotopic maps of simplicial sets are   homotopic morphisms of cdgl's.
\end{enumerate}
\end{corollary}


\section{Free extensions and cofibrations}
Here we present a large class of cofibrations.

\begin{definition} A  \emph{free extension} of a cdgl $L$ is an inclusion
$$L\hookrightarrow L'  = (L\, \widehat{\amalg}\, \widehat{\mathbb L}(V),d)$$ for which the following properties are satisfied:  \begin{enumerate}

\item[$\bullet$] $V = V_{\geq -1} $ and $V_{-1}$ is generated by Maurer-Cartan elements.
\item[$\bullet$] $V_0= V_0' \oplus V_0''$, where each element of $V_0'$ is a cycle and $V_0''$ has a basis $\{x_i\}_{i\in I}$ for which there exist Maurer-Cartan elements $a_i$ and $b_i$ in $L_{-1}$ or $V_{-1}$ such that $(\widehat{\mathbb L}(a_i, b_i, x_i),d)$ is a Lawrence-Sullivan interval included in $L'$.
\item[$\bullet$] For $x\in V_n$, with $n\geq 1$, there is a Maurer-Cartan element $a$ such that $d_ax\in L\,\widehat{\amalg}\,{\widehat{\mathbb L}}(V_{<n})$.
\end{enumerate}
\end{definition}

\begin{theorem}\label{free}
Every free extension is a cofibration.
\end{theorem}

Note that this class of cofibrations contains the classical Koszul-Quillen extension
\cite[Section II.5]{tan} in the simply connected setting.

\begin{proof} We need to prove that for every commutative square
$$\xymatrix{
 L\ar[r]\ar[d]
 &
 A\ar[d]^-{p}\\
 L'\ar[r]_{\varphi}
 \ar@{-->}[ru]^\phi
 &
  B
  }$$
in which $p$ is both a fibration and a weak equivalence, there exists the dotted morphism $\phi$ making commutative both triangles.

Recall that the cdgl morphisms
$$
\lasu_{\mathscr I}=\{\lasu_{\partial\Delta^n}\hookrightarrow \lasu_{\Delta^n}\}_{n\ge0},\quad
\lasu_{\mathscr J}=\{\lasu_{\wedge^n_i}{\hookrightarrow} \lasu_{\Delta^n}\}_{n\ge0, \,i=0,\dots,n},
$$
are generating sets of cofibrations and trivial cofibrations, respectively. 
In particular $L\hookrightarrow (L\, \widehat{\amalg}\, \widehat{\mathbb L}(V_{-1}\oplus V_0''),d)$ is a cofibration. 
Indeed, this morphism is obtained by successive pushouts of diagrams $L\leftarrow 0\to\lasu_{\Delta^0}$ 
(one for each generator of $V_{-1}$) and 
$L\leftarrow \lasu_{\partial\Delta^1}\to\lasu_{\Delta^1}$ (one for each $i\in I$), 
and by definition, this is a generated cofibration. 
Thus, there exists a cdgl morphism $\phi\colon  (L\, \widehat{\amalg}\, \widehat{\mathbb L}(V_{-1}\oplus V_0''),d)\to A$ 
as in the diagram.

On the other hand, since $d(V_0')=0$, $p$ is surjective in non-negative degrees, and $p^0\colon A^0\to B^0$ is a 
quasi-isomorphism, $\phi$ is easily extended to
$(L\, \widehat{\amalg}\, \widehat{\mathbb L}(V_{-1}\oplus V_0),d)$.

We finish by defining $\phi$ inductively on $V_{\ge 1}$. Assume $\phi$ defined on $V_{<n}$ with $n\geq 1$ and 
let $x\in V_{n}$. We denote by $a$ the  Maurer-Cartan 
element for which $d_ax\in L\,\widehat{\amalg}\,{\widehat{\mathbb L}}(V_{<n})$. 
In the restriction  to the component of $a$,
$$\xymatrix{
& A^{\phi(a)}\ar[d]^p\\
\bigl(L\,\widehat{\amalg}\,{\widehat{\mathbb L}}(V_{<n})\bigr)^a  \ar[r]^-{\varphi}& B^{\varphi(a)},
}
$$
the element $\phi d_{a} ({x})$ is a cycle in $A^{\phi(a)}$ with $p(\phi d_{a}({x})) = d_{\varphi(a)}\varphi(x)$.
Therefore, since $p$ is surjective in non-negative degrees and the restriction $p\colon A^{\phi(a)}\to B^{\varphi(a)}$ is a  quasi-isomorphism,   
there exists $y\in A$ with $\phi d_{a}({x})= d_{\phi(a)}y$ and $p(y)= \varphi({x})$. 
We define $\phi({x}) = y$ and  observe that
$\phi(dx)=\phi(d_{a}x)-\phi([a,x])=d_{\phi(a)}y- [\phi(a),y]=dy=d\phi(x)$.
\end{proof}

\begin{definition}\label{def:freecofibrant}
A free extension of~0 is called a \emph{free cofibrant} cdgl.
\end{definition}

 \begin{corollary}
 The model $\lasu_X$ of any finite simplicial set is a free cofibrant\break  cdgl.
 \end{corollary}

We end this section with examples of weak equivalences.
 Let $U$ be a   graded vector space and let $sU$ be the suspension of $U$,
  $|su|=|u|+1$.
 We form the   cdgl $(\libc(U\oplus sU),d)$ with $du=0$ and $dsu=u$, for any $u\in U$. Then,

\begin{proposition}\label{prop:we2}
For any cdgl $A$, the  projection on the first factor,
$$p\colon  A\, \widehat{\amalg}\, (\libc(U\oplus sU),d)\to A,$$
is a trivial fibration.
\end{proposition}

\begin{proof} 
The map $p$ being surjective is a fibration and we are reduced to prove that it is a weak equivalence.

Write $B=A\, \widehat{\amalg}\, (\libc(U\oplus sU),d)$.
Denote by $\cF_{p}A$ the complete descending filtration of $A$. We endow $\hL(U\oplus sU)$
with its filtration as completion of a free Lie algebra, i.e. by the closure of the descending central series. We
define now the complete descending filtration of $B$, $\cF_{p}B$, as coproduct of cdgl's.

On the other hand, we denote by $B(n)$ the sub vector space of $B$
generated by the brackets containing exactly $n$ 
elements in $U\oplus sU$.  Each $B(n)$ is a sub complex and ker$\, p$ is the product
 $\prod_{n=1}^\infty B(n)$.

We define a derivation $s$ of degree $+1$ in $B$   by $s(A)= 0$, $s(u) = su$ and $s(su)= 0$.
The derivation $\theta=sd+ds$ is equal to the identity on $U\oplus sU$ and so is  the multiplication by $n$ on $B(n)$.
 
 The kernel $K_p$ of $F_p(B)\to F_p(A)$ is also a product of sub complexes,
  $\prod_{n=1}^{\infty} B(n)\cap K_p$. 
 Therefore if $\alpha = \sum_{n} \alpha_n$ is a cycle with $\alpha_n \in B(n)\cap K_p$, 
 then $$\alpha_n = \frac{1}{n} (ds\alpha_n + sd\alpha_n) = d(\frac{1}{n} s\alpha_n)$$
and $$\alpha = d(\sum_n \frac{1}{n} s\alpha_n)\,.$$  
From Theorem~\ref{propDR} it follows that $p$ is a weak equivalence. 
 \end{proof}

\begin{corollary}\label{cor:we2}
 The canonical injection,
$\iota\colon A\to A\, \widehat{\amalg}\, (\libc(U\oplus sU),d)$,
is a weak equivalence.
\end{corollary}

\begin{proof}
The equality $p\circ \iota=\id$ and the fact that $p$ is a weak equivalence imply that $\iota$ is a weak equivalence.
 \end{proof}

\section{The cylinder construction and the Lawrence-Sullivan interval}

We define a cylinder construction for cdgl's along the same line as the original definition for dgl's   in 
\cite[Section II.5]{tan}.

\medskip
 Consider a cdgl of the form $(\libc(V),d)$. We denote by $U$ a copy of $V$,
 the isomorphism  $V\xrightarrow{\cong} U$ being represented by $v_i\mapsto u_i$, where $v_i$ and $u_i$ are respectively a graded basis for $V$ and for $U$.
 We construct the cdgl,
 $$(\libc(V\oplus U\oplus sU),d),$$
 where $d_{|_V}$ is the differential on $(\libc(V),d)$ and, for $u\in U$,   $du=0$ and $dsu=u$.

 A derivation $i$ of degree $+1$ is defined on $(\libc(V\oplus U\oplus sU),d)$
 by $i(v) = su$, $i(u)= i(su)= 0$. Then $\theta = i\circ d + d\circ i$ is a derivation commuting with $d$.
 Therefore $e^{\theta}$ is an automorphism of $(\libc(V\oplus U\oplus sU),d)$. We
define graded vector spaces $V'$ and $\overline{V}$ respectively isomorphic to $V$ and $sV$, and define a  morphism of graded Lie algebras,
 \begin{equation}\label{equa:lepsi}
 \psi \colon
 \libc(V\oplus  {V}'\oplus \ov{V})\to \libc(V\oplus U\oplus sU),
 \end{equation}
 by $\psi(v)=v$, $\psi(v')=e^{\theta}(v)$ and $\psi(\ov{v})=su$.
 This  is an isomorphism that induces a differential
 $D = \psi^{-1}\circ d\circ \psi$ on $\libc(V\oplus V'\oplus \ov{V})$. Since $e^\theta$ is an automorphism commuting with the differential $d$, the sub Lie algebra
$ \libc(V')$ is a sub  cdgl, isomorphic to $(\libc(V),d)$.

\begin{definition}\label{def:cylinder}
Let $L$ be a cDGL of the form $L=(\libc(V),d)$.
The \emph{cylinder construction} on $L$ is the cDGL,
$$\cyl(L)=  (\libc(V\oplus V'\oplus \ov{V}),D)\,$$
together with the maps
$$\xymatrix@1{
L
\ar@<2pt>[r]^-{\iota_{0}}
\ar@<-2pt>[r]_-{\iota_{1}}
&
\cyl(L)\ar[r]^-{p}
&
L,
}$$
defined by
$\iota_{0}(v)=v$, $\iota_{1}(v)=v'$, $p(v)=p(v')=v$, $p(\ov{v})=0$.
\end{definition}

\begin{proposition} The cylinder {\em Cyl}$(L)$ is a cylinder object in the model category {\bf cDGL}.  \end{proposition}

\begin{proof} From the isomorphism (\ref{equa:lepsi}), we deduce
$$\mbox{Cyl}(L) \cong L \,\widehat{\amalg}\,\, \widehat{\mathbb L} (U\oplus sU).$$
The result follows from Proposition~\ref{prop:we2} and Theorem~\ref{free}.
 \end{proof}

By applying Theorem~\ref{LST}, we find the first example of a cylinder construction.

\begin{corollary}
The Lawrence-Sullivan interval $\lasu_{\Delta^1}$ is isomorphic to  the cylinder on $\lasu_{\Delta^0}$,
$$\lasu_{\Delta^1} \cong \cil\, \lasu_{\Delta^0}\,.$$
\end{corollary}

Since the cylinder constructed above is a cylinder in the model category {\bf cDGL}, we have

\begin{corollary} Two morphisms $f, g \colon (\widehat{\mathbb L}(V),d) \to L'$ of
cdgl's are \emph{homotopic} 
if there exists a morphism $F \colon\cil\,(\widehat{\mathbb L}(V),d) \to L'$ such that 
$f= F\circ \iota_0$ and $g= F\circ \iota_1$. 
\end{corollary}

 \vspace{3mm} In a  cdgl  $(\widehat{\mathbb L}(V),d)$ we write $d = \sum_{i\geq 1} d_i$ where $d_i \colon V\to \mathbb L^i(V)$. If $f \colon (\widehat{\mathbb L}(V),d)\to (\widehat{\mathbb L}(W,d)$ is a morphism of dgl's, we write $f = \sum_{i\geq 1}f_i$ where $f_i(V)\subset \mathbb L^i(W)$.

\begin{proposition} Let $f,g \colon (\widehat{\mathbb L}(V),d)\to (\widehat{\mathbb L}(W),d)$ be two homotopic maps then the induced maps
 $$f_1, g_1 \colon (V, d_1) \to (W, d_1)$$
 are homotopic maps of chain complexes.
 \end{proposition}

 \begin{proof} Let $H\colon \cil(\widehat{\mathbb L}(V),d) \to ({\widehat{\mathbb L}}(W),d)$ be a homotopy between $f$ and $g$.  Since $D_1(\overline{v}) = v'-v - \overline{d_1v}$, we have
 $$g_1v - f_1v = H_1 \overline{d_1v} + d_1H_1\overline{v}\,.$$
 The morphism $h \colon V\to W$ defined by $h(v) = H_1\overline{v}$ is then a chain homotopy between $f_1$ and $g_1$. \end{proof}

\begin{proposition} Two   homotopic maps $f,g\colon L\to L'$  induce the same map on the set of 
equivalence classes of Maurer-Cartan elements.\end{proposition}

\begin{proof} By definition, there is a morphism $H \colon \mbox{Cyl}(L)\to L'$ such that $f = H\circ \iota_0$ and $g = H\circ \iota_1$. Now let $a$ be a MC element in $L$. The naturalness of the cylinder construction implies the commutativity of the diagram
$$\xymatrix{
{\lasu}_0 \ar[d]_a \ar@/^/[rr]^{i_a}\ar@/_/[rr]_{i_b} && {\lasu}_{\Delta^1}\ar[d]^{\cyl (a)}\\
L \ar@/^/[rr]^{\iota_0}\ar@/_/[rr]_{\iota_1} && \cyl (L)\ar[r]^(0.58)H &L'.}$$
Here $i_a$ and $i_b$ denote the injections of the Maurer-Cartan elements $a$ and $b$ of ${\lasu}_{\Delta^1}$. 
Then $H\circ \mbox{Cyl}(a)$ is a path between $f(a)$ and $g(a)$.
\end{proof}

\section{Minimal Models of connected finite simplicial sets}\label{sec:modelsset}

In this section, we introduce the notion of minimal model of a connected finite simplicial set, establish its unicity up to isomorphism 
 and present some examples.
 
 \begin{definition}[\emph{Cofibrant model}]
Let $L$ be a cdgl, and $\{a_i\}$ be a set of Maurer-Cartan elements representing the different 
gauge equivalent classes. 
Since $L^{a_i}$ is concentrated in degrees $\geq 0$, we can construct a graded vector space 
$V(i)$ with $V(i) = V(i)_{\geq 0}$ and a quasi-isomorphism
$$\varphi_i : (\widehat{\mathbb L}(V_i), d_{a_i}) \longrightarrow L^{a_i}\,.$$
For degree reasons, $d(V(i)_n) \subset \widehat{\mathbb L}(V(i)_{<n})$.
The union of the $\varphi_i$ induces a morphism
$$\varphi \colon L'= (\widehat{\mathbb L}((\oplus_i \Q a_i)\oplus (\oplus_i V(i)), d) \to L$$
where $d(a_i)= -\frac{1}{2}[a_i, a_i]$ and for $x\in V(i)$, $dx= d_{a_i}x - [a_i,x]$. By construction $\varphi$ is a weak equivalence and, by Theorem \ref{free}, $L'$ a cofibrant replacement, or {\em model}, of $L$.\end{definition}

\begin{definition}\label{def:minimal}
A cdgl $(\widehat{\mathbb L}(V),d)$ is \emph{minimal} if $V = V_{\geq 0}$ and $d(V)\subset \mathbb L^{\geq 2}(V)$.
\end{definition}

 If $(L,d)$ is a cdgl such that $H(L,d)=H_{\geq 0}(L,d)$, then, there 
  exists a minimal cdgl $(\widehat{\mathbb L}(V),d)$ equipped with a quasi-isomorphism
$$\varphi \colon (\widehat{\mathbb L}(V),d) \to (L,d)\,.$$
This minimal cdgl is called a \emph{Lie minimal model} of $(L,d)$.

\medskip
For the geometrical application we are developping in this section, 
we focus on Lie minimal models of cofibrant cdgl's. 
In particular, we determine the vector space of generators of these models.

\begin{theorem}\label{thm:minimodel}
Any homologically connected cdgl $(\hL(Z),d)$ is isomorphic to a coproduct
$$(\hL(Z),d)\cong (\hL(V),d) \widehat{\amalg}\; \hL(E\oplus dE),$$
where $(\hL(V),d)$ is a minimal model of $(\hL(Z),d)$ and $V\cong H_{*}(Z,d_{1})$, with $d_{1}$
the linear part of the differential $d$.
Moreover, two Lie minimal models of $(\hL(Z),d)$ are isomorphic.
\end{theorem}

\begin{proof}
We decompose the differential $d$ along the bracket length as $d=\sum_{i\geq 1}d_{i}$, with $d_{i}(Z)\subset \L^i(Z)$.
As $d_{1}$ is a differential, we can write
$Z=V\oplus E\oplus d_{1}E$ with $(d_{1})_{|V}=0$. Therefore, the projection
$$p\colon (\hL(Z),d_{1})\to (\hL(V),0),$$
defined by $p_{|E}=p_{|d_{1}E}=0$, is a quasi-isomorphism. Its kernel $K$ is thus $d_{1}$-acyclic.
We define a grading on it by $K^i=K\cap \L^i(Z)$.

We use an induction on this grading and suppose that the differential $d$ verifies
$d(V)\subset \hL(V)\oplus K^{\geq n}$. Let $v\in V$, we decompose 
$$d_{n+1}(v)=\alpha(v)+\beta(v)\in \L^{n+1}(V)\oplus K^{n+1}.$$
As $d^2=0$ and $(d_{1})_{|V}=0$, the   induction implies
$$d_{1}d_{n+1}v=-\sum_{i=2}^nd_{i}d_{n+2-i}v\in \L^{n+1}(V).$$
Thus, the component $\beta(v)$ is a $d_{1}$-cycle and there exists $\gamma(v)\in K^{n+1}$ such that
$\beta(v)=d_{1}\gamma(v)$. By doing this argument on
each element of a basis $(v_{i})_{i}$ of $V$, we construct a graded
vector space $V'$, isomorphic to $V$, generated by $v_{i}\mapsto v_{i}-\gamma(v_{i})$. By construction, $V'$
verifies $d(V')\subset \hL(V')\oplus K^{\geq n+1}$.
An iteration gives a graded vector space, that we still denote $V$, such that 
$dV\subset \hL(V)$.

To conclude, we observe that the injection
$ \hL(V\oplus E\oplus dE)\to \hL(Z)$
has a linear part which is an isomorphism, and so is an isomorphism. 
The first part of the statement is now a consequence of Corollary~\ref{cor:we2}.

\medskip
For the second property, let
$(\hL(W),d)\to (\hL(Z),d)$ be a minimal model. Its composition with the previous projection $p$ gives
a weak-equivalence $f\colon (\hL(W),d)\to (\hL(V),d)$. Recall  $W=W_{\geq 0}$,  $V=V_{\geq 0}$ and
observe that it is sufficient to prove that the linear part, $f_{1}$, of $f$ is an isomorphism.
We use an induction on the degrees in $V$ and $W$. Suppose that $f_{1}\colon W_{<n}\to V_{<n}$ is an isomorphism.

First, let $v\in V_{n}$.  As $dv\in \hL(V)_{n-1}$ is a boundary, 
there exists $u\in \hL(W)_{n-1}$ such that $f(du)=dv$. 
The element $v-f(u)$ being a cycle in $\hL(V)$, there exist $w\in\hL(W)$ and $x\in \hL(V)$ such that
$$f(w)=v-f(u)+dx.$$
This equality implies $f_{1}(u+w)=v$ and the surjectivity of $f_{1}\colon W_{n}\to V_{n}$.

To establish the injectivity, let $w\in W_{n}$ such that $f_{1}(w)=0$. As $f_{1}$ is surjective, we can find 
$w'\in \L^{\geq 2}(W)$ such that $f(w+w')=0$. Thus $f(d(w+w'))=0$ and, by induction, $d(w+w')=0$. The map $f$
being a quasi-isomorphism, there exists $x\in W$ such that 
$$w+w'=dx.$$
As $d_{1}=0$ and $w'$ is decomposable, we deduce $w=0$, as expected.
\end{proof}

\begin{definition}\label{def:minimalss}
Let $X$ be a connected simplicial set (or a path connected topological space).
 Denote by $a$ a base point in $X$. 
 The Lie minimal model, 
 $\varphi_{X}\colon m_{X}=(\widehat{\mathbb L}(V),d)\to (\lasu_X^a,d_a)$ is called \emph{the Lie minimal model of $X$.}
 By construction $V\cong s^{-1} H_*(X;\mathbb Q)$. 
\end{definition}

Moreover, if $f\colon X\to Y$ is a continuous map between path connected spaces, then $f$ has a model 
$m_f\colon m_X\to m_Y$ and we have a diagram of Lie minimal models that commutes up to homotopy 
$$\xymatrix{ 
(\lasu_X^a, d_a) \ar[rr]^{\lasu_f} && (\lasu_Y^{f(a)}, d_{f(a)})\\
m_X = (\widehat{\mathbb L}(V),d) \ar[u]^\varphi \ar[rr]^{m_f} && m_Y= (\widehat{\mathbb L}(W),d). \ar[u]^{\psi}
}$$
By construction the linear part of $m_f$ is the desuspension of $H_*(f;\mathbb Q)$.

\begin{example}[\emph{Wedge of circles}]\label{exam:wedgecircle}
Let $X=\vee_{i=1}^n S^1$ be a wedge of circles. From $\ov{H}_{*}(X;\Q)\cong H_{1}(X;\Q)$
and Theorem~\ref{thm:minimodel}, the minimal model of $X$ is
$$(\hL(\oplus_{i=1}^n \Q a_{i}),0)$$
with $|a_{i}|=0$.
\end{example}

Our machinery allows also the computation of the Lie minimal model of   a space $X$
obtained by attaching a 2-cell along a wedge $W$ of  circles.
In fact $X$ is the union of a wedge of $p$ circles $S^1_{(1)}, \dots , S^1_{(p)}$ with a $n$-gone,
in which each face is identified with one of the circles.
Let
 $\omega = x_{i_1}^{\varepsilon_1} \dots x_{i_n}^{\varepsilon_n}$ be the word associated to the 
 adding of the 2-cell, 
 with $\varepsilon_i = \pm 1$ in function of the orientations.
The fundamental group of $X$  is given by the presentation
$$\pi_1(X) = \langle \, x_1, \dots , x_p ; \omega\, \rangle,$$
and its minimal model has a similar form.

\begin{proposition}\label{prop:2cell}
With the previous notation, the  Lie minimal model of  $X$ 
is the cdgl $(\hL(x_1, \dots , x_p, y),d)$ with $|x_{i}|=0$,
$|y|=1$,
$d x_{i}=0$ and $d  y = x_{i_1}^{\varepsilon_1} * \dots * x_{i_n}^{\varepsilon_n}$,
where  $*$ denotes the Baker-Campbell-Hausdorff product.
\end{proposition}

For instance, the Lie minimal model of  the Klein Bottle is
$  (\hL(u,v,y),d)$
with $| u| = | v| = 0$, $| y| = 1$, $du= d v= 0$  and $$d y= u*v*u*v^{-1}.$$

\begin{proof}
With Theorem~\ref{thm:minimodel}, the Lie minimal model of the $n$-gone, whose boundary is a circuit of $n$ circles,
 is defined on the complete Lie algebra
 $\hL(b_{1},\dots,b_{n},z)$
 with $|b_{i}|=0$ and $|z|=1$. From Example~\ref{exam:wedgecircle}, we know  that $db_{i}=0$
and \cite[Theorem 2]{lawsu} implies that $dz= b_1*\dots * b_n$.

On the other hand, the Lie minimal model of $X$ is defined on the complete Lie algebra,
$(\hL(x_{1},\dots,x_{p},y)$ with $|x_{i}|=0$ and $|y|=1$. As before, we know that $dx_{i}=0$.
 
  We also have a natural map $f \colon Z\to X$ mapping each circle of the boundary of $Z$ on one of the circles
  in the wedge $W$. 
  At the level of the induced map $m_{f}$ between the Lie minimal models, this means
  $m_{f}(b_{i})=x_{i_{j}}^{\varepsilon_{i}}$.
   Moreover $H_{2}(f;\Q)$ is an isomorphism. 
Thus there is a nonzero rational number $\lambda$ such that $y':= m_f(z)= \lambda y+ \omega$, where $\omega$ is a decomposable element.
Replacing $y$ by $y'$ gives an isomorphism
$m_{X}\cong (\widehat{\mathbb L}(x_1, \dots , x_p, y'),d)$.
Since $m_f$ commutes with the differential,  we have
$dy' = m_fdz  = x_{i_1}^{\varepsilon_1} * \dots * x_{i_n}^{\varepsilon_n}$.
 \end{proof}

\begin{corollary}
Let $G= <x_1, \dots , x_p \colon \omega>$ be a finitely generated one-relation group, with 
$\omega = x_{i_1}^{\varepsilon_1} *\dots * x_{i_n}^{\varepsilon_n}$ and $\varepsilon_{i}=\pm 1$. 
Then the rational Malcev completion of $G$ is the group
$$\widehat{\mathbb L} (a_1, \dots , a_p)/ ( a_{i_1}^{\varepsilon_1} *\dots *a_{i_p}^{\varepsilon_p})\,,$$
where $|a_{i}|=0$ and the group law is the Baker-Campbell-Hausdorff product.
\end{corollary}

\begin{proof} This follows directly from Proposition~\ref{prop:2cell} and the fact that 
$H_0(\lasu_X)$ is the Malcev completion of $\pi_1(X)$, see \cite[Theorem 9.1]{pri}.
\end{proof}

\section{Relation with other model category structures  on {\bf cDGL}}\label{sec:fdgl}

Here, we compare our model structure on $\catcdgl$ with other known model structures.

\medskip
First, one may consider in $\catcdgl$ the classical model structure given  on categories of unbounded chain complexes
 enriched with some algebraic structure, see for instance \cite[\S2]{hi}. 
 Fibrations are surjective morphisms, weak equivalences are quasi-isomorphisms and cofibrations 
 are morphisms satisfying the left lifting property with respect to trivial fibrations.

The zero map $0\to \lib(a)$, in which $a$ is a Maurer-Cartan element, is not surjective but 
it is a fibration in our model structure. The same example is a quasi-isomorphism but it is not a weak equivalence in our structure.

\medskip
On the other hand, in \cite{lamar}, A. Lazarev and M. Markl    define  a  model category structure on 
the full subcategory of $\catcdgl$
formed by the  profinite complete dgl's.

\begin{definition}
 A dgl $L$ is \emph{profinite-complete} if $L= \varprojlim_\alpha L_\alpha$, where the $L_\alpha$'s are 
 nilpotent finite dimensional dgl's. 
 We denote by $\catcdgl^f$ the associated subcategory of $\catcdgl$.
The \emph{profinite-completion} of a dgl $L$ is the projective limit 
$$\widehat{L}^f= \varprojlim_\alpha L_\alpha$$
where $L_\alpha$ runs over all the nilpotent finite dimensional quotients of $L$.
\end{definition}

 \begin{example} Let $L$ be an infinite dimensional abelian Lie algebra $L = \oplus_{n}\, \mathbb Q e_n$ 
 concentrated in degree $k\geq 0$. 
 Then $L$ is complete, but not profinite-complete.  Its profinite-completion is the product $\prod_n \mathbb Q e_n$. 
 \end{example}

\begin{example} 
 1) Let $V$ be a graded vector space and $\{e_\alpha\}$ be a  basis for $L=\mathbb L(V)$,
 as graded vector space.
  Then $\widehat{\mathbb L}^f(V):= \widehat{L}^f $ is isomorphic, as  graded vector space, to 
  $  \prod_\alpha \mathbb Q e_\alpha$. 
  When $V$ is finite dimensional, then $\widehat{\mathbb L}(V)\cong \widehat{\mathbb L}^f(V)$. 
 
 \medskip
 2) Denote by $L$ the Lie algebra concentrated in degree 0, obtained as the quotient of the free Lie 
 algebra $\mathbb L(\oplus_{i\geq 1}\Q a_i\oplus_{i\geq 1} \Q b_i)$ by the ideal  $I$ 
  generated by the brackets $[a_i, a_j]$, $[b_i, b_j]$ and $[a_i, b_j]$ for $i\neq j$ and by the elements
   $[a_i,b_i]-[a_1,b_1]$ for $i\geq 2$. 
   Denote by $\varphi\colon L\to E$  a morphism from $L$ in a finite dimensional Lie algebra $E$ and let 
   $a_n + \sum_{i<n}\lambda_i a_i+\sum_j\mu_{j}b_{j}$ be an element in $\ker \varphi$. Then
 $$0 = \varphi\left[a_n+ \sum_{i<n}\lambda_i a_i, b_n\right]= 
 \varphi[a_n, b_n]\,.$$
 As $\widehat{L}^f$ is built as a projective limit of finite dimensional Lie algebras, the canonical map
  $p\colon L\to \widehat{L}^f$ is not injective and $[a_1, b_1] \in \ker p$. 
\end{example}

Remark 7.2 of \cite{lamar} proves the following.

\begin{lemma}
When $L$ is a finitely generated dgl, then $\widehat{L}^f = \widehat{L}$.
\end{lemma}

We defined a functor $\coch^*\colon \catcdgl^f \to \bf{CDGA}$ 
as a generalization of the usual cochain functor, $C^*$, from the category 
{\bf DGL} of connected dgl's to the category  of augmented commutative differential graded algebras (cdga's).
 More precisely, let $L= \varprojlim_{\alpha} L(\alpha)$ be a profinite-complete dgl. 
 Since  each $L(\alpha)$ is a finite dimensional  nilpotent dgl, 
 we can form the usual cochain algebra $C^*(L(\alpha))$ and set  
$$\coch^*(L) = \varinjlim_{\alpha} C^*(L(\alpha))\,.$$
Remark also that if $L$ is the dual of some dglc $E$, then $\coch^*(L)= {\mathcal A}(E)$. 

We consider also a functor ${\mathscr L}^*\colon \bf{CDGA} \to \catcdgl^f$ as the composition of the functor
${\mathcal E}\colon \catcdga\to {\bf DGLC}$ previously introduced with a  duality: if $A$ is a cdga, one sets
$${\mathscr L}^* (A)=({\mathcal E}(A))^{\sharp}$$

\medskip
The functors ${\mathscr L}^*$ and $\coch^*$ have
been defined   in \cite{lamar}, where the authors  prove the following  (Propositions 7.5 and 9.10 in op. cit.).

\begin{theorem}\label{markl} With the above notations,
\begin{enumerate}
\item[(i)] The functors $\coch^*$ and ${\mathscr L}^*$ are adjoint.
\item[(ii)] For any augmented cdga $A$, the counit $\coch^*{\mathscr L}^* (A)\to A$ is a quasi-isomorphism,
\item[(iii)] Let $L$ be a profinite-complete dgl and $\varphi_L {\colon} {\mathscr L}^* \coch^*(L) \to L$ the associated counit, then $\coch^*(\varphi_L) {\colon} \coch^*(L) \to \coch^*{\mathscr L}^*\coch^*(L)$ is a quasi-isomorphism.
\end{enumerate}\end{theorem}

A model category structure on $\catcdgl^f$ is defined  in \cite{lamar} by:
 \begin{itemize}
\item a morphism $f$ is a \emph{weak equivalence} if $\coch^*(f)$ is a quasi-isomorphism,
\item $f$ is a \emph{fibration} if it is a surjection,
\item $f$ is a \emph{cofibration} if it has the left lifting property with respect to all trivial fibrations.
\end{itemize}

The class of fibrations with this model structure is properly contained in the class of our fibrations.
 The class of weak equivalences is different  from ours. 
For instance let $L$ be the abelian Lie algebra on one generator $x$ in degree $-3$, and let 
$f\colon 0\to L$ be the injection. Clearly $f$ is a weak equivalence
in the sense of Theorem~\ref{model2}, 
because $\widetilde\mc(L)=\{0\}$ and $L^0=0$. However $\coch^*(f)$ is not a quasi-isomorphism.

The main result of this section is a comparison theorem for the two notions of weak equivalences:
they coincide for a large class of maps, including the models of simplicial maps.
\emph{In what follows the terms of cofibration, fibration and weak equivalence will always refer 
to the model category defined in Section~\ref{sec:closedmodel}.}

\begin{theorem}\label{propcompar} 
Let $f {\colon} (L,d)=(\widehat{\mathbb L}(V),d)\to (L',d)= (\widehat{\mathbb L}(V'),d)$ be a morphism of free cofibrant  cdgl's with $V$ and $V'$ finite dimensional. Then  $\coch^*(f)$ is a quasi-isomorphism if and only if  $f$ is a weak equivalence.
\end{theorem}

For the proof we  need some preliminaries on $\mathbb Z$-graded Sullivan algebras, and  properties 
of   ${\mathscr L}^*$.

\begin{definition} A \emph{($\mathbb Z$-graded) Sullivan algebra} is a cdga of the form $(\Lambda V,d)$, where  
$\Lambda V$ is the free graded commutative algebra on a $\mathbb Z$-graded vector space $V$
 endowed with an increasing filtration, $V(0)\subset V(1)\subset \dots $, $V = \cup_{n\geq 0} V(n)$  
  such that $dV(0)= 0$ and for $n>0$,
$$dV(n)\subset \Lambda V(n-1)\,.$$
\end{definition}

For any profinite-complete dgl $L$, the cochain algebra $\coch^*(L)$ is a $\mathbb Z$-graded Sullivan algebra.

\begin{definition}\label{def:homotopycdga}
 Two morphisms $h,k \colon (\Lambda V,d)\to (\Lambda W,d)$ between 
 $\mathbb Z$-graded Sullivan algebras are \emph{homotopic} if there is a morphism,
$$H {\colon} (\Lambda V\otimes \Lambda \overline{V}\otimes \Lambda \widehat{V}, D)\to (\Lambda W,d),$$
where $D\overline{v} = \widehat{v}$, $D\widehat{v}= 0$, $h(v)= H(v)$ and
$$k(v) = H(e^{sD+Ds})(v)\,.$$
Here, $s$ is the derivation defined by $sv= \overline{v}$, $s\overline{v}= s\widehat{v} = 0$.
\end{definition}

The cdga $(\Lambda V\otimes \Lambda \overline{V}\otimes \Lambda \widehat{V},D)$ is a cylinder object. We have in particular a diagram
$$\xymatrix{(\Lambda V,d) \ar@<2pt>[r]^-{\iota_{0}}
\ar@<-2pt>[r]_-{\iota_{1}}
& (\Lambda V\otimes \Lambda \overline{V}\otimes \Lambda\widehat{V},D) \ar[r]^(0,65)p & (\Lambda V,d)}$$
such that $p\circ i_0=p\circ i_1= \id$, $i_0(v)= v$, $i_1(v)= e^{sD+Ds}(v)$, $p(v)= v$, $p(\overline{v})=p(\widehat{v})= 0$.

\medskip
 Consider now the effect of ${\mathscr L}^*$ on the cylinder object of the $\mathbb Z$-graded cdga $(\Lambda V,d)$.
Since the ideal $I$
  generated by $\overline{V}$ and $\widehat{V}$ is acyclic, as a graded vector space, 
  we can  write $I = S\oplus D(S)$ where $D {\colon} S\to D(S)$ is an isomorphism.
Therefore we obtain an isomorphism
$$ {\mathscr L}^*(\Lambda V\otimes \Lambda \overline{V}\otimes \Lambda \widehat{V},D)  \stackrel{\cong}{\longrightarrow}  {\mathscr L}^*(\Lambda V,d)\,\widehat{\amalg}^f\, \libc^f(E\oplus dE)\,,$$
where $E = D(S)^{\sharp}$ and $\widehat{\amalg}^f$ defines the coproduct in the category of profinite-complete dgl's. 

We need the analog of Proposition~\ref{prop:we2} in the profinite context.

\begin{proposition}\label{prop:weprofinite}
For any profinite-complete dgl $M$, the  projection on the first factor,
$$p\colon  M\, \widehat{\amalg}\, (\libc^f(E\oplus dE),d)\to M,$$
is a trivial fibration.
\end{proposition}

\begin{proof}
Set $L=M\, \widehat{\amalg}\, (\libc^f(E\oplus dE),d)$ and 
denote by $B_n$ the sub complex of $L$ generated by the Lie brackets containing exactly $n$ terms from $E\oplus dE$. 
We define a derivation, $s$ of $L$ by $s(M)=s(E)= 0$ and $s(de)=e$ for $e\in E$. 
Then $sd+ds$ is the multiplication by $n$ in $B_n$.
 
The kernel $K$ of $p$  is the product, $K= \prod_{i\geq 1} B_i$. Therefore an element $a\in K$ can be written as
$a= \sum_{n\geq 1} a_n$ with $a_n\in B_n$. 
If $a$ is a cycle then we have $a= d(-\sum_{n\geq 1} \frac{sa_n}{n})$ and $a$ is a boundary. 
 We have proved that  $p$ is a quasi-isomorphism.
 
 Now, we   consider the Maurer-Cartan elements.
Let $\alpha$ be a Maurer-Cartan element in $L$. We write it as $\alpha = \sum_{n\geq 0} \alpha_n$ with 
$\alpha_n \in B_n$. 
We first observe that $\alpha_{0}$ is a Maurer-Cartan element in $M$.

We use  an induction on $n$, supposing that there are elements $b_i\in B_i$, for $1\leq i\leq n$, such that 
$$(b_n*b_{n-1}*\dots *b_1)\cG \alpha -\alpha_0\in \prod_{q>n}B_q.$$ 
Set $b = b_n*b_{n-1}* \dots * b_1$, and write
 $$b\cG \alpha = \alpha_0 + \beta_{n+1}+ \beta_{n+2}+ \dots$$
 with $\beta_i \in B_i$. Since $b\cG\alpha$ is a Maurer-Cartan element, we have  
 $d\beta_{n+1}= -[\beta_{n+1}, \alpha_0]$. 
 Applying the derivation $s$, we get, 
 $$\beta_{n+1}= \frac{1}{n+1} ( ds\beta_{n+1} - [s\beta_{n+1}, \alpha_0]).$$
 Set $b_{n+1}= s\beta_{n+1}$. 
 Then a  computation from Definition~\ref{def:BCH}  shows that 
 $$(b_{n+1}*b)\cG \alpha -\alpha_0\in \prod_{q>n+1} B_q\,.$$
 The sequence $b_1$, $b_2*b_1, \dots $ converges in $L$ to an element $x$  
 because it converges in all the nilpotent quotients $L_\alpha$. 
 Moreover $x\cG \alpha = \alpha_0$ as this is true in all the quotients $L_\alpha$. 
 
 We have proved that the natural injection $M\to (L,d)$ induces a surjection 
 $\widetilde{\MC}(M)\to \widetilde{\MC}(L,d)$. Since  $M$ is a retract of $(L,d)$ the corresponding map 
 $$\widetilde{\MC}(M) \to \widetilde{\MC}(L,d)$$ is a bijection.
\end{proof}

By the same argument as in Corollary \ref{cor:we2},  
the morphisms ${\mathscr L}^*(p)$,  ${\mathscr L}^*(i_0)$ and ${\mathscr L}^*(i_1)$
are weak equivalences.

\begin{corollary}\label{corhom}
Let $h,k{\colon} (\Lambda V,d)\to (\Lambda W,d)$ be homotopic maps between $\mathbb Z$-graded Sullivan algebras. Then $\widetilde\mc({\mathscr L}^*(h)) = \widetilde\mc({\mathscr L}^*(k))$ and for each Maurer-Cartan element $a$ in ${\mathscr L}^*(\Lambda W,d)$, $H({\mathscr L}^*(h), d_a)= H({\mathscr L}^*(k), d_a)$.
\end{corollary}

For the proof of Theorem~\ref{propcompar}, we need the following result.

\begin{proposition} \label{prop:LCqi}
 For any free cofibrant cdgl $(L,d)= (\widehat{\mathbb L}(V),d)$ with $\dim V<\infty$,  the counit $\varphi {\colon} {\mathscr L}^*\coch^*(L) \to L$ is a weak equivalence.
\end{proposition}  

\begin{proof}  We equip $L$ with the filtration $L = F_1 \supset F_2\supset \cdots$ where $F_p$ is the subspace generated by the iterated brackets of length $\geq p$. This is a complete descending filtration. 
We denote by $G(L)$ the graded Lie algebra associated to the filtration $F_p$, $G_p= F_p/F_{p+1}$. Then, 
 $L= \oplus G_p$ and $F_p = \oplus_{q\geq p} G_q$.
  By construction $d(G_p)\subset \oplus_{q\geq p} G_q$ and $[G_p, G_q]\subset G_{p+q}$.

Let us write $(\Lambda V,d) =\coch^*(L,d)$, where $V = sL^{\sharp}$. We equip   $V$ with the induced grading, 
$V_p= s(G_p^{\sharp})$ and extend it in an algebra grading to $\Lambda V$.  
By construction $d(\Lambda V)_q \subset (\Lambda V)_{\leq q}$.     

Apply now the functor ${\mathscr L}^*$,  $${\mathscr L}^*(\Lambda V,d) = (\libc^f(X),d)$$
 with $X  = s^{-1}(\Lambda V)^{ \sharp}$.  This induces a  grading on $X$ by $X_p = s^{-1}((\Lambda V)_p)^{\sharp}$
 that we extend in a Lie algebra grading on $\hL^f(X)$.   
 The associated filtration
 $$F_p({\mathscr L}^*(\Lambda V)) = \oplus_{q\geq p} ({\mathscr L}^*(\Lambda V))_q$$ 
 satisfies $d(F_p)\subset F_p$ and $[F_p, F_q]\subset F_{p+q}$. 
 Moreover, this is a complete descending filtration,
  $${\mathscr L}^*(\Lambda V) = \varprojlim_n {\mathscr L}^*(\Lambda V)/ F_n.$$ 
The counit
$\varphi {\colon} {\mathscr L}^*\coch^*(L)\to L$ is a dgl morphism defined on the generators by 
$\varphi(s^{-1}(\Lambda^{\geq 2}V)^{\sharp}) = 0$ and $\varphi(s^{-1}(sx^{\sharp})^{\sharp}) = x$. 
Therefore, the kernel $K$ of $\varphi$ is  generated by $s^{-1}(\Lambda^{\geq 2}V)^{\sharp}$
and by the elements
$$[s^{-1}(sa^{\sharp})^{\sharp}, s^{-1}(sb^{\sharp})^{\sharp}]- s^{-1}(s[a,b]^{\sharp})^{\sharp}$$
for $a,b\in  L$.
Clearly $\varphi$ is a morphism of  filtered  dgl's.

 In $L$ and  ${\mathscr L}^*\coch^*(L)$, the differential $d$ can be written as $d = d_{0}+d_{1}+ \dots$,
  where $d_i$ increases the new grading by $i$.   
  In $(L,d_0)$ the differential is linear and we can modify the degrees in order that all the elements are in positive degree. 
  Denote by $L'$ this new cdgl. 
  With this transformation ${\mathscr L}^*\coch^* (L)$ becomes $\mathcal L^*C^*(L')$, 
  where $\mathcal L^*$ and $C^*$ 
  are the usual Quillen functors. 
  Since  
  $\varphi \colon \mathcal L^*C^*(L')\to L'$
  is a quasi-isomorphism, the same is true for $(L,d_0)$.  

Now by induction on $p$, the five lemma argument shows that the induced map
$$\varphi_p : {\mathscr L}^*\coch^*(L)/F_p\to L/F_p$$ is a quasi-isomorphism. Therefore for each integer $p$, the induced map $F_p({\mathscr L}^*\coch^*(L))\to F_p(L)$ is a quasi-isomorphism. 
 By  Theorem \ref{propDR}, $\widetilde{\mc}(\varphi)$ is a bijection.
 
 The same argument now shows that for any Maurer-Cartan element in ${\mathscr L}^*\coch^*(L)$, 
 the induced map $\varphi \colon ({\mathscr L}^*\coch^*(L),d_a)\to (L, d_{\varphi (a)})$ is a quasi-isomorphism. 
 It follows that $\varphi$ is a weak equivalence. 
 \end{proof}

\begin{example}
 Let $(L,d)= (\mathbb L(a),d)$ where $a$ is a Maurer-Cartan element.  
 Write $x= sa^{\sharp}$ and $y= s[a,a]^{\sharp}$. 
 In $(\Lambda V,d)= \coch^*(L,d)$,  we have
$$dy= x-x^2\,.$$
Now form ${\mathscr L}^*(\Lambda V,d)= (\widehat{\mathbb L}^f(X),d)$ and denote by  $a_n$ and $b_n$ the elements in $X$ defined by  $a_n = s^{-1}(x^n)^{\sharp}$ and $b_n= s^{-1}(yx^n)^{\sharp}$. In $(\libc^f(X),d)$ we have 
$$da_1 = b_1\,, da_2 = b_2-b_1-\frac{1}{2} [a_1, a_1]\,,$$
$$da_n = b_n-b_{n-1} = \frac{1}{2} \sum_{1\leq i\leq n-i\leq n-1} [a_i, a_{n-i}].$$
Denoting $\omega = \sum_{i=1}^\infty a_i$,  we have
$$d\omega = -\frac{1}{2}[\omega, \omega] \hspace{5mm}\mbox{and } \varphi (\omega) = a\,.$$
 This example shows that the profinite-completion was absolutely necessary in the construction of the functor 
 ${\mathscr L}^*$ for having a surjection between the Maurer-Cartan sets, $\widetilde{\MC}(-)$.
\end{example}

\begin{proof}[Proof of Theorem~\ref{propcompar}]
 By construction $\coch^*(L)$ and $\coch^*(L')$ are Sullivan algebras. Now by the same argument as in 
 \cite[Proposition 14.6]{FHTI} $\coch^*(f)$ admits an inverse $g$ up to homotopy. 
 Therefore, by Corollary \ref{corhom}, $\widetilde{\MC}({\mathscr L}^*\coch^*f)$ is a bijection 
 with inverse $\widetilde{\MC}({\mathscr L}^*(g))$ and for each MC element, 
 ${\mathscr L}^*\coch^*(f)$ is a quasi-isomorphism. 
 This means that ${\mathscr L}^*\coch^*(f)$ is a weak equivalence.

Now, the commutativity of the diagram
$$\xymatrix{{\mathscr L}^*\coch^*L \ar[r]^(0.60)\varphi\ar[d]_{{\mathscr L}^*\coch^*f} & L\ar[d]^f\\
{\mathscr L}^*\coch^*L' \ar[r]^(0.60)\varphi & L'}$$ and the fact that $\varphi$ is a weak equivalence imply that $f$ is a weak equivalence.

To prove the converse, we first study the kernel  $K $  of the projection 
$(\widehat{\mathbb L}^f(V\oplus U\oplus sU),d)\to (\widehat{\mathbb L}^f(V),d)$. 
The derivation $s$ on $\widehat{\mathbb L}^f(V\oplus U\oplus sU)$ defined by $s (u)=su, s (v)= 0$ and 
$s(su)=0$ satisfies $s d+ds= 0$ on $K$ and on each quotient 
$K/K\cap \widehat{\mathbb L}^{\geq p}(V\oplus U\oplus sU)$. 
Therefore $H^*(\coch^*(K))= 0$. Remark now that 
the cylinder of a cdgl $L$ has a image by $\coch^*$ that can be decomposed as
$$\coch^*(\mbox{Cyl}(L)) \cong   \coch^*(\widehat{\mathbb L}^f(V\oplus U\oplus sU),d))
\cong \coch^*(L) \otimes \coch^*(K)\,.$$
Therefore $\coch^*(p)$ is a quasi-isomorphism, and since $p\circ i_0=p\circ i_1= \id$, 
it follows that $\coch^*(i_0)=\coch^*(i_1)$ induce the same map in homology. 
As a consequence,  if
two maps $k, \ell \colon (L,d) \to (L',d)$ are homotopic, there is a map $H\colon \mbox{Cyl}(L)\to L'$ with
 $\ell = H\circ i_0$ and $k= H\circ  i_1$
 and we deduce $H^*(\coch^*(\ell)) = H^*(\coch^*(k))$.

Let $f$ be a weak equivalence between cofibrant objects. As usual in a model category, 
$f$ admits an inverse $g$ up to homotopy. 
Since $f\circ g$ and $g\circ f$ are homotopic to the identity, 
$H^*(\coch^*(f))$ and $H^*(\coch^*(g))$ are inverse of each other. 
In particular, $\coch^*(f)$ is a quasi-isomorphism.
\end{proof}

\providecommand{\bysame}{\leavevmode\hbox to3em{\hrulefill}\thinspace}
\providecommand{\MR}{\relax\ifhmode\unskip\space\fi MR }
\providecommand{\MRhref}[2]{%
  \href{http://www.ams.org/mathscinet-getitem?mr=#1}{#2}
}
\providecommand{\href}[2]{#2}

\end{document}